\documentclass[preprint,12pt,numbers]{elsarticle}

\usepackage[a4paper,margin=1in]{geometry}
\usepackage{amsmath,amssymb,amsthm,amsfonts}
\usepackage{mathtools}
\usepackage{bm}
\usepackage{enumitem}
\usepackage{graphicx}
\usepackage{hyperref}
\usepackage{booktabs}
\usepackage{color}
\usepackage{float}
\usepackage{cases}
\usepackage{subcaption}

\numberwithin{equation}{section}

\newtheorem{theorem}{Theorem}[section]

\newtheorem{remark}{Remark}[section]

\newcommand{\uh}{\bm u}
\newcommand{\uw}{\bm u_w}
\newcommand{\utau}{\bm u_\tau}
\newcommand{\us}{\bm u_s}
\newcommand{\vh}{\bm v}
\newcommand{\nh}{\bm n}
\newcommand{\gradtau}{\nabla_\tau}
\newcommand{\dd}{\,\mathrm d}
\newcommand{\uhtau}{\bm u_{h,\tau}}
\newcommand{\vhtau}{\bm v_{h,\tau}}

\newcommand{\tiluhtaus}{\bm {\tilde{u}}_{h,\tau}^{*}}
\newcommand{\utildehtaup}{\widetilde{\bm u}_{h,\tau}^{n+1}}

\journal{Communications in Nonlinear Science and Numerical Simulation}
\begin{document}

\begin{frontmatter}



\title{A Fully Decoupled Unified SAV–ZEC Finite Element Scheme for the Phase-Field Moving-Contact-Line Model}


\author[aff1]{Xiangjun Meng}
\ead{xjmeng@bjut.edu.cn}

\author[aff2]{Hanlin Guo}
\ead{hanlin_guo@mail.bnu.edu.cn}

\author[aff3]{Xuelian Bao\corref{cor1}}
\ead{baoxuelian@scut.edu.cn}

\cortext[cor1]{Corresponding author.}
\affiliation[aff1]{
	organization={School of Mathematics, Statistics and Mechanics, Beijing University of Technology},
	city={Beijing},
	postcode={100124},
	country={China}
}

\affiliation[aff2]{
	organization={School of Mathematical Sciences, Beijing Normal University},
	city={Beijing},
	postcode={100875},
	country={China}
}

\affiliation[aff3]{
	organization={School of Mathematics, South China University of Technology},
	city={Guangzhou},
	postcode={510641},
	state={Guangdong},
	country={China}
}

\begin{abstract}
We propose a fully decoupled, linear, second-order and energy-stable finite element scheme for the Cahn--Hilliard--Navier--Stokes system with the generalized Navier boundary condition and dynamic contact-line condition, which couples phase transport, capillary forces, fluid convection, and contact-line motion. Specifically, we introduce a unified scalar auxiliary variable (SAV) for the sum of the bulk and boundary energies and a zero-energy-contribution (ZEC) variable for the nonlinear coupling terms. 
Combining this reformulation with BDF2 time discretization and a pressure-correction procedure yields a scheme that requires only solving linear elliptic problems and updating scalar auxiliary variables at each time step. For stationary walls, we establish an unconditional discrete energy inequality for a modified energy
that is bounded from below. Numerical experiments confirm the second-order temporal accuracy and the discrete energy law. Several two- and three-dimensional simulations further illustrate the performance of the method for moving-contact-line problems on both rectangular and curved computational domains.
\end{abstract}

%


\end{frontmatter}




\section{Introduction}
 Moving contact line problems involve coupling between two-phase flow, interfacial forces, and fluid–solid interactions \cite{bonn2009,snoeijer2013}. When a fluid--fluid interface moves along a solid wall, contact-line motion depends on capillary forces, wall slip, wall friction, and contact-angle relaxation. In the classical sharp-interface formulation, the no-slip condition leads to the moving-contact-line singularity \cite{huh1971}, while phase-field models avoid this difficulty by replacing the sharp interface with a diffuse transition layer and by incorporating wall wettability through a boundary free energy \cite{anderson1998,jacqmin2000,ding2007,yue2010,qian2003,qian2006}. 
 The resulting Cahn--Hilliard--Navier--Stokes (CH-NS) system, supplemented with the generalized Navier boundary condition and the dynamic contact-line condition, provides a thermodynamically consistent model for moving contact lines \cite{qian2006,shenyang2010,yu2017,ren2007}.

A main difficulty in constructing an energy-stable scheme for the system is the treatment of the nonlinear coupling terms in the bulk and on the boundary. The velocity transports the phase field, while the chemical potential contributes a capillary force to the momentum equation. On the solid wall, the dynamic contact-line condition and the generalized Navier boundary condition couple the phase field and the tangential velocity. 

Existing methods address different parts of this difficulty. Convex-splitting, stabilization, and fractional-step formulations can preserve an energy structure, but many require coupled or nonlinear solves \cite{eyre1998,shen2010,gao2012,he2011,salgado2013}. Linear energy-stable moving-contact-line schemes based on explicit--implicit splitting were developed in \cite{shen2015,alandchen2016,ma2017}, while continuous finite element discretizations were studied in \cite{jiang2015}. In particular, second-order linear stabilized discretizations have been developed for Cahn--Hilliard models with dynamic boundary conditions, together with rigorous stability and convergence analyses \cite{meng2023ch,meng2025}. IEQ and SAV methods provide systematic frameworks for linearizing nonlinear free energies \cite{yang2016,shen2018}; their moving-contact-line extensions include second-order IEQ schemes \cite{yangyu2018}. 
Fully decoupled multiple-SAV schemes have also been developed for related CH--NS systems \cite{li2022}.
Nevertheless, in moving-contact-line models, the nonlinear energy contains both bulk and boundary contributions, and the hydrodynamic coupling contains both interior and boundary terms.  
Recently, Ref. \cite{gao2025} introduced separate SAVs for the bulk and boundary energies and a ZEC variable to derive a linear, fully decoupled BDF2 scheme, which requires splitting only with respect to the ZEC variable. They established a time-discrete energy inequality and used staggered-grid finite differences in space. 

The objective of this work is to construct a second-order finite element scheme for the moving-contact-line model that is linear, fully decoupled, and easy to implement. To this end, we introduce a unified SAV formulation that uses a single scalar auxiliary variable to represent the combined bulk and boundary energies. This formulation enables the nonlinear double-well potential in the domain and the nonlinear boundary energy on the solid boundary to be treated simultaneously.  It also yields quadratic
auxiliary-variable terms in the modified discrete energy,
which has an explicit lower bound.
We also introduce a zero-energy-contribution (ZEC) variable whose evolution equation contains the transport, capillary, convection, and contact-line coupling terms. These terms cancel in the continuous energy estimate. The auxiliary equation is discretized so that the corresponding cancellations also hold at the discrete level.

Based on this SAV-ZEC reformulation, we construct a BDF2 finite element scheme combined with a pressure-correction method for the incompressibility constraint \cite{chorin1967,guermond2006,he2022}.  The nonlinear coefficients are extrapolated explicitly, while the phase-field variables, auxiliary variables, intermediate velocity, and pressure are determined sequentially. The resulting scheme is linear and fully decoupled, and its discrete energy structure is inherited from the reformulated system.

The proposed method differs from previous approaches in the following two respects. First, the auxiliary energy variable is not assigned only to the bulk potential; instead, it is constructed for the combined bulk and boundary free energy, which is particularly suitable for dynamic wetting problems. 
Notice that, we combine the nonlinear bulk and boundary energies into a single energy functional and introduce one scalar auxiliary variable, rather than treating them separately. This unified treatment reduces the number of auxiliary variables and evolution equations, simplifying the SAV–ZEC reformulation and its numerical implementation while preserving energy stability and a linear, decoupled solution structure. Moreover, it only requires that the combined nonlinear energy be bounded from below, without imposing separate lower-bound assumptions on the two contributions.

Second, the finite element formulation accommodates boundary-fitted triangulations of non-rectangular domains and incorporates the contact-line conditions in weak form, allowing simulations on complex geometries. And the ZEC treatment of bulk and boundary coupling terms decouples the phase-field and velocity subproblems while preserving a discrete energy law. 
Numerical experiments are designed to verify the convergence rate, energy dissipation, and mass conservation. In addition to standard two- and three-dimensional shear-driven droplet simulations	\cite{shen2015,ma2017,yu2017,yangyu2018}, we also include relaxation tests on curved solid boundaries \cite{bao2012} to demonstrate the robustness of the method in geometrically nontrivial moving contact line configurations. 

The main contributions of this work are summarized as follows.
\begin{enumerate}[label=(\roman*)]
    \item The nonlinear bulk and boundary energies are treated through
    a unified SAV, reducing the number of auxiliary
    	evolution equations and yielding a quadratic auxiliary energy. The ZEC equation accounts for the coupling terms in both
    	the bulk equations and the contact-line conditions.
    	
    \item The resulting BDF2 scheme decouples the phase-field and velocity computations. Within each time step, the subproblems share the same coefficient matrix, while the auxiliary variables are obtained by solving scalar linear equations.

    \item For stationary walls, we establish an unconditional discrete energy dissipation law for a modified energy bounded from below. The dissipation terms account for phase diffusion, viscous dissipation, contact-line relaxation, and Navier slip.

    \item Convergence and conservation tests are supplemented by
    shear-driven droplet simulations in two and three dimensions.
    A relaxation test on a polygonal approximation of a curved
    wall further illustrates the application of the finite element
    scheme to non-rectangular domains.
\end{enumerate}
The rest of this paper is organized as follows. In Section~\ref{sec:model}, we formulate the CH--NS moving-contact-line model and recall its energy-dissipation law. The SAV-ZEC reformulation is derived in Section~\ref{sec:reformulation}. Based on this reformulation, the fully discrete BDF2 finite element scheme is developed in Section~\ref{sec:scheme} and its unconditional energy stability is established in Section~\ref{sec:stability}. The corresponding decoupled solution procedure is detailed in Section~\ref{sec:implementation}. Numerical results are reported in Section~\ref{sec:numerics}, followed by concluding remarks in Section~\ref{sec:conclusion}.

\section{The phase-field moving contact line model}
\label{sec:model}
We first recall the phase-field moving-contact-line model and derive its energy law, which is the basis for the numerical scheme developed below.
Consider a bounded polygonal or polyhedral domain $\Omega\subset\mathbb R^d$ with $d=2$ or $3$. 
We denote by $\Gamma\subset\partial\Omega$ the solid wall on which contact-line motion may occur. 
For clarity of presentation, we assume that the two fluids have equal densities and viscosities.
The unknowns are the phase-field variable $\phi$, the chemical potential $\mu$, the fluid velocity $\uh$, and the pressure $p$. 
In the bulk domain, the evolution is governed by the CH-NS system
\begin{subequations}
	\label{eq:original}
	\begin{numcases}{}
		\phi_t+\nabla\cdot(\uh\phi)=M\Delta\mu, \label{eq:ch}\\
		\mu=\lambda\bigl(-\varepsilon\Delta\phi+f(\phi)\bigr), \label{eq:mu}\\
		\uh_t+(\uh\cdot\nabla)\uh+\nabla p-\nu\Delta\uh+\phi\nabla\mu=0, \label{eq:ns}\\
		\nabla\cdot\uh=0. \label{eq:div}
	\end{numcases}
\end{subequations}
Here $M>0$ is the mobility, $\nu>0$ is the viscosity, $\varepsilon>0$ denotes the interfacial thickness, and $\lambda>0$ is the capillary coefficient. 
The bulk free energy density is given by the classical double-well potential
\begin{equation*}
	F(\phi)=\frac{1}{4\varepsilon}(\phi^2-1)^2,
	\qquad
	f(\phi)=F'(\phi)=\frac{1}{\varepsilon}(\phi^3-\phi).
\end{equation*}
Thus, the Cahn--Hilliard equation describes the diffusive relaxation of the phase interface, while the Navier--Stokes equation is driven by the capillary force induced by the chemical potential.

The interaction between the diffuse interface and the solid wall is described through a boundary energy
\begin{equation*}
	g(\phi)=-\frac{\sqrt2}{3}\cos\theta_s
	\sin\left(\frac{\pi}{2}\phi\right),
\end{equation*}
where $\theta_s$ is the prescribed static contact angle. 
This boundary energy gives rise to the boundary functional
\begin{equation*}
	L(\phi)=\varepsilon\partial_n\phi+g'(\phi),
\end{equation*}
which measures the deviation of the dynamic contact angle from its equilibrium state.

On the solid wall $\Gamma$, the system is supplemented by the generalized Navier boundary condition (GNBC) and the dynamic contact-line condition (DCLC):
\begin{subequations}
	\label{eq:bc}
	\begin{numcases}{}
		\uh\cdot\nh=0, 
		& on $\Gamma$, \label{eq:un}\\
		\nu l(\phi)(\utau-\uw)
		+\nu\partial_n\utau
		-\lambda L(\phi)\gradtau\phi=0,
		& on $\Gamma$, \label{eq:gnbc}\\
		\partial_n\mu=0,
		& on $\Gamma$, \label{eq:mun}\\
		\phi_t+\utau\cdot\gradtau\phi=-\gamma L(\phi),
		& on $\Gamma$, \label{eq:dclc}\\
		L(\phi)=\varepsilon\partial_n\phi+g'(\phi),
		& on $\Gamma$. \label{eq:L}
	\end{numcases}
\end{subequations}
Here $\nh$ is the unit outward normal vector, $\utau=\uh-(\uh\cdot\nh)\nh$
is the tangential velocity, and
$
\gradtau v=\nabla v-(\partial_n v)\nh
$
denotes the tangential gradient. 
The nonnegative function $l(\phi)$ denotes the wall-friction coefficient, $\uw$ is the prescribed tangential wall velocity satisfying $\uw\cdot\nh=0$ on $\Gamma$, and $\gamma>0$ is the contact-line relaxation parameter. Moreover, the impermeability condition $\uh\cdot\nh=0$ implies $\uh=\utau$ on $\Gamma$.

A key property of the above model is that the bulk evolution and the wall dynamics satisfy a common energy law.
The total energy consists of the kinetic energy, the diffuse-interface energy in the bulk, and the boundary energy:
\begin{equation*}
	E_{\rm tot}(\uh,\phi)
	=
	\frac12\|\uh\|^2
	+\lambda\left(
	\frac{\varepsilon}{2}\|\nabla\phi\|^2
	+\int_\Omega F(\phi)\dd x
	+\int_\Gamma g(\phi)\dd s
	\right).
\end{equation*}
To derive the energy law, we take the $L^2$-inner products of \eqref{eq:ch}, \eqref{eq:mu}, and \eqref{eq:ns}
with $\mu$, $-\phi_t$, and $\uh$, respectively. 
After integration by parts, the boundary contributions from the bulk equations are exactly balanced by the GNBC and DCLC terms in \eqref{eq:bc}. 
This yields the energy identity
\begin{equation}
	\label{eq:continuous-energy}
	\frac{\dd}{\dd t}E_{\rm tot}
	=
	-\nu\|\nabla\uh\|^2
	-M\|\nabla\mu\|^2
	-\lambda\gamma\|{L(\phi)}\|_\Gamma^2
	-\nu\|{\sqrt{l(\phi)}\,\us}\|_\Gamma^2
	-\nu(l(\phi)\us,\uw)_\Gamma,
\end{equation}
where $\us=\utau-\uw$ denotes the slip velocity. 
The first four terms on the right-hand side correspond to viscous dissipation, Cahn--Hilliard diffusion, contact-line relaxation, and Navier-slip dissipation, respectively. 
The remaining term accounts for the mechanical work input from the moving wall.

\begin{remark}
    If $\uw=0$, the last term in \eqref{eq:continuous-energy}
    vanishes, yielding an energy dissipation law.
\end{remark}

\section{The unified SAV-ZEC reformulation}
\label{sec:reformulation}

We reformulate the model by introducing a unified SAV for the
nonlinear bulk and boundary energies and a ZEC variable for the
nonlinear coupling terms. With consistent initial data, the
resulting system is equivalent to the original model and provides
the basis for a linear, decoupled scheme satisfying a discrete
energy law.

The first component is the unified SAV formulation for the nonlinear
bulk and boundary energy. We combine these contributions into a
single energy functional,
\begin{equation*}
	\mathcal E_1(\phi)
	:=
	\int_\Omega F(\phi)\dd x
	+\int_\Gamma g(\phi)\dd s
	+A,
\end{equation*}
where the constant $A>0$ is chosen such that $\mathcal E_1(\phi)>0$. This requires only a lower bound on the combined nonlinear
bulk and boundary energy, allowing greater flexibility than separate lower-bound assumptions on the two contributions.
The corresponding scalar auxiliary variable is defined by
\begin{equation*}
	U(t)=\sqrt{\mathcal E_1(\phi)}.
\end{equation*}
For later use, we introduce
\begin{equation*}
	H(\phi)
	=
	\frac{f(\phi)}
	{\sqrt{\mathcal E_1(\phi)}},
	\qquad
	J(\phi)
	=
	\frac{g'(\phi)}
	{\sqrt{\mathcal E_1(\phi)}}.
\end{equation*}
By differentiating $U(t)$ with respect to time, one obtains
\begin{equation*}
	U_t
	=
	\frac12(H(\phi),\phi_t)
	+
	\frac12(J(\phi),\phi_t)_\Gamma.
\end{equation*}
Hence, one scalar variable is sufficient to represent the nonlinear bulk--boundary energy contribution. 
This unified treatment is particularly effective for moving-contact-line model because the boundary energy directly determines the contact-angle relaxation through
\begin{equation}
	L(\phi)=\varepsilon\partial_n\phi+g'(\phi),\notag
\end{equation}
rather than serving merely as an auxiliary boundary correction.

The second component is a zero-energy-contribution variable. 
The purpose of this variable is not to alter the continuous model, but to account for the nonlinear coupling terms that cancel in the energy estimate.
We introduce $Q(t)$ by
\begin{equation}
	\label{eq:Qode}
	\begin{aligned}
		Q_t={}&
		\int_\Omega
		\left[
		\nabla\cdot(\uh\phi)\mu
		+\phi\nabla\mu\cdot\uh
		+((\uh\cdot\nabla)\uh)\cdot\uh
		\right]\dd x        \\
		&-\lambda(L(\phi)\gradtau\phi,\utau)_\Gamma
		+\lambda(L(\phi)\gradtau\phi,\utau)_\Gamma,
		\qquad
		\mbox{with } Q(0)=1.
	\end{aligned}
\end{equation}
At the continuous level, the right-hand side of \eqref{eq:Qode} has zero contribution to the energy variation: the transport and capillary terms cancel with each other, the hydrodynamic convection term vanishes in the energy estimate, and the two boundary contact-line terms cancel pairwise. 
Thus the exact solution satisfies $Q(t)=1$. 
Retaining $Q$ as an auxiliary
unknown permits explicit treatment of the nonlinear coupling
terms while preserving the cancellations in the discrete
energy estimate.

With the unified SAV variable $U$ and the ZEC variable $Q$, the original system can be rewritten as follows:
\begin{subequations}
	\label{eq:reform}
	\begin{numcases}{}
		\phi_t+Q\nabla\cdot(\uh\phi)=M\Delta\mu,
		\label{eq:reform-ch}\\
		\mu=\lambda\left(
		-\varepsilon\Delta\phi
		+H(\phi)U
		\right),
		\label{eq:reform-mu}\\
		U_t=\frac12(H(\phi),\phi_t)
		+\frac12(J(\phi),\phi_t)_\Gamma,
		\label{eq:reform-U}\\
		\uh_t+Q(\uh\cdot\nabla)\uh+\nabla p-\nu\Delta\uh
		+ Q\phi\nabla\mu=0,
		\label{eq:reform-u}\\
		\nabla\cdot\uh=0,
		\label{eq:reform-div}\\
		\begin{aligned}
			Q_t={}&
			\int_\Omega
			\left[
			\nabla\cdot(\uh\phi)\mu
			+\phi\nabla\mu\cdot\uh
			+((\uh\cdot\nabla)\uh)\cdot\uh
			\right]\dd x\\
			&-\lambda(L(\phi)\gradtau\phi,\utau)_\Gamma
			+\lambda(L(\phi)\gradtau\phi,\utau)_\Gamma.
		\end{aligned}
		\label{eq:reform-Q}
	\end{numcases}
\end{subequations}
The corresponding boundary conditions are
\begin{subequations}
	\label{eq:reform-bc}
	\begin{numcases}{}
		\phi_t+Q\utau\cdot\gradtau\phi=-\gamma L(\phi),
		& on $\Gamma$, \label{eq:reform-bc-dclc}\\
		L(\phi)=\varepsilon\partial_n\phi+J(\phi)U,
		& on $\Gamma$, \label{eq:reform-bc-L}\\
		\uh\cdot\nh=0,
		& on $\Gamma$, \label{eq:reform-bc-un}\\
		\nu l(\phi)(\utau-\uw)
		+\nu\partial_n\utau
		-\lambda QL(\phi)\gradtau\phi=0,
		& on $\Gamma$, \label{eq:reform-bc-gnbc}\\
		\partial_n\mu=0,
		& on $\Gamma$. \label{eq:reform-bc-mu}
	\end{numcases}
\end{subequations}

This reformulated system is equivalent to the original model at the continuous level. 
Indeed, the consistent initial condition gives $U^2=\mathcal E_1(\phi)$, while the zero-energy-contribution equation gives $Q(t)=1$. 
Substituting these identities into \eqref{eq:reform}--\eqref{eq:reform-bc} recovers \eqref{eq:original}--\eqref{eq:bc}. 

\begin{theorem}[Energy law of the reformulated system]
	Assume that $\uw=0$. 
	Then the reformulated system \eqref{eq:reform}--\eqref{eq:reform-bc} satisfies
	\begin{equation}
		\label{eq:reform-energy-law}
		\frac{\dd}{\dd t}\mathcal E(\uh,\phi,U,Q)
		=
		-\nu\|\nabla\uh\|^2
		- M\|\nabla\mu\|^2
		-\lambda\gamma\|L(\phi)\|_\Gamma^2
		-\nu\|\sqrt{l(\phi)}\utau\|_\Gamma^2,
	\end{equation}
	where
	\begin{equation*}
		\label{eq:reform-energy}
		\mathcal E(\uh,\phi,U,Q)
		=
		\frac12\|\uh\|^2
		+\lambda\left(
		\frac{\varepsilon}{2}\|\nabla\phi\|^2
		+|U|^2
		\right)
		+\frac12|Q|^2-\lambda A.
	\end{equation*}
\end{theorem}

\begin{proof}
	Taking the $L^2(\Omega)$ inner products of \eqref{eq:reform-ch}, \eqref{eq:reform-mu}, and \eqref{eq:reform-u} with $\mu$, $-\phi_t$, and $\uh$, respectively, gives the bulk part of the energy balance. 
	The equation for $U$ is multiplied by $2\lambda U$, which converts the nonlinear bulk and boundary free-energy contributions into the time derivative of $\lambda |U|^2$. 
	Similarly, multiplying the equation for $Q$ by $Q$ gives the terms needed to cancel the nonlinear coupling terms in the energy estimate.
	
	For the boundary terms, the dynamic contact-line condition gives
	\begin{align}
		\lambda\varepsilon(\partial_n\phi,\phi_t)_\Gamma
		&=
		\lambda(L(\phi)-J(\phi)U,\phi_t)_\Gamma
		\nonumber\\
		&=
		\lambda(L(\phi),
		-Q\utau\cdot\gradtau\phi-\gamma L(\phi))_\Gamma
		-\lambda U(J(\phi),\phi_t)_\Gamma
		\nonumber\\
		&=
		-\lambda Q(L(\phi)\gradtau\phi,\utau)_\Gamma
		-\lambda\gamma\|L(\phi)\|_\Gamma^2
		-\lambda U(J(\phi),\phi_t)_\Gamma.
		\label{eq:proof-boundary-phi}
	\end{align}
	On the other hand, the generalized Navier boundary condition yields, for $\uw=0$,
	
	\begin{equation}\label{eq:proof-boundary-u}
		\nu(\partial_n\utau,\utau)_\Gamma
		=
		\lambda Q(L(\phi)\gradtau\phi,\utau)_\Gamma
		-\nu\|\sqrt{l(\phi)}\utau\|_\Gamma^2.
	\end{equation}
	The first term on the right-hand side of \eqref{eq:proof-boundary-phi} is exactly canceled by the corresponding term in \eqref{eq:proof-boundary-u}. 
	The term involving $U(J(\phi),\phi_t)_\Gamma$ is canceled by the boundary part of the $U$ equation. 
	The transport, capillary, convection, and boundary coupling terms are canceled by the ZEC equation for $Q$. 
	After these cancellations, the only remaining terms are the dissipative contributions in \eqref{eq:reform-energy-law}. 
	This completes the proof.
\end{proof}

\section{Finite element weak formulation and fully discrete scheme}
\label{sec:scheme}

We next introduce the finite element setting and the fully discrete scheme. 
Let $\mathcal T_h$ be a conforming and shape-regular triangulation or tetrahedralization of $\Omega$, consisting of open disjoint elements $K$ such that
\[
\overline\Omega=\bigcup_{K\in\mathcal T_h}\overline K .
\]
For an integer $l\ge0$, $P_l(K)$ denotes the space of polynomials on $K$ with total degree no larger than $l$.

We use the following conforming finite element spaces:
\begin{equation*}
	\label{eq:fe-spaces}
	\begin{aligned}
		Y_h
		&=
		\left\{
		\varphi_h\in C^0(\overline\Omega):
		\varphi_h|_K\in P_{l_1}(K),\ \forall K\in\mathcal T_h
		\right\},\\
		V_h
		&=
		\left\{
		\vh_h\in C^0(\overline\Omega)^d:
		\vh_h|_K\in P_{l_2}(K)^d,\ \forall K\in\mathcal T_h
		\right\}
		\cap H^1_{\nh}(\Omega)^d,\\
		O_h
		&=
		\left\{
		q_h\in C^0(\overline\Omega):
		q_h|_K\in P_{l_2-1}(K),\ \forall K\in\mathcal T_h
		\right\}
		\cap L_0^2(\Omega),
	\end{aligned}
\end{equation*}
where
\begin{equation*}
	H^1_{\nh}(\Omega)^d
	=
	\left\{
	\vh\in H^1(\Omega)^d:
	\vh\cdot\nh=0 \ \text{on } \Gamma
	\right\},
	\qquad
	L_0^2(\Omega)
	=
	\left\{
	q\in L^2(\Omega):
	\int_\Omega q\dd x=0
	\right\}.
\end{equation*}
Thus, $Y_h\subset H^1(\Omega)$, $V_h\subset H^1_{\nh}(\Omega)^d$, and $O_h\subset L_0^2(\Omega)$. 
The auxiliary variables $U(t)$ and $Q(t)$ are nonlocal scalar variables and are therefore treated as real-valued unknowns rather than finite element functions.

For later use, we denote by $Y_h|_\Gamma$ the trace space of $Y_h$ on $\Gamma$, and define the tangential trace space
\[
V_{h,\tau}|_\Gamma
=
\left\{
\vhtau=\vh_h-(\vh_h\cdot\nh)\nh:
\vh_h\in V_h
\right\}.
\]
Before giving the fully discrete scheme, we derive the finite element weak formulation of the reformulated system. 
Find
\[
\phi_h,\mu_h\in Y_h,\qquad
\uh_h\in V_h,\qquad
p_h\in O_h,\qquad
U,Q\in\mathbb R,
\]
such that for all $\omega_h,\eta_h\in Y_h$, $\vh_h\in V_h$, and $q_h\in O_h$,
\begin{subequations}
	\label{eq:semi-weak-reform}
	\begin{numcases}{}
		((\phi_h)_t,\omega_h)
		-Q(\uh_h\phi_h,\nabla\omega_h)
		=
		-M(\nabla\mu_h,\nabla\omega_h),
		\label{eq:semi-weak-phi}\notag\\[1mm]
		\begin{aligned}
			(\mu_h,\eta_h)
			=
			\lambda\Big[
			&\varepsilon(\nabla\phi_h,\nabla\eta_h)
			-\varepsilon(\partial_n\phi_h,\eta_h)_\Gamma
			+U(H(\phi_h),\eta_h)
			\Big],
		\end{aligned}
		\label{eq:semi-weak-mu}\notag\\[1mm]
		U_t
		=
		\frac12(H(\phi_h),(\phi_h)_t)
		+\frac12(J(\phi_h),(\phi_h)_t)_\Gamma,
		\label{eq:semi-weak-U}\notag\\[1mm]
		\begin{aligned}
			((\uh_h)_t,\vh_h)
			&+Q((\uh_h\cdot\nabla)\uh_h,\vh_h)
			+\nu(\nabla\uh_h,\nabla\vh_h)
			-(p_h,\nabla\cdot\vh_h)
			-\nu(\partial_n\uhtau,\vhtau)_\Gamma\\
			&+Q(\phi_h\nabla\mu_h,\vh_h)
			=0,
		\end{aligned}
		\label{eq:semi-weak-u}\notag\\[1mm]
		(\nabla\cdot\uh_h,q_h)=0,
		\label{eq:semi-weak-div}\notag\\[1mm]
		\begin{aligned}
			Q_t={}&
			-((\uh_h\phi_h),\nabla\mu_h)
			+(\phi_h\nabla\mu_h,\uh_h)
			+((\uh_h\cdot\nabla)\uh_h,\uh_h)\\
			&-\lambda(L(\phi_h)\gradtau\phi_h,\uhtau)_\Gamma
			+\lambda(L(\phi_h)\gradtau\phi_h,\uhtau).\notag
		\end{aligned}
		\label{eq:semi-weak-Q}
	\end{numcases}
\end{subequations}
Here
\[
\vhtau=\vh_h-(\vh_h\cdot\nh)\nh,
\qquad
\uhtau=\uh_h-(\uh_h\cdot\nh)\nh .
\]
The boundary conditions on $\Gamma$ are imposed in the following weak form:
\begin{subequations}
	\label{eq:semi-weak-bc}
	\begin{numcases}{}
		((\phi_h)_t,m_h)_\Gamma
		+Q(\uhtau\cdot\gradtau\phi_h,m_h)_\Gamma
		=
		-\gamma(L(\phi_h),m_h)_\Gamma,
		& for all $m_h\in Y_h|_\Gamma$, \label{eq:semi-weak-bc-dclc}\notag\\[1mm]
		(L(\phi_h),r_h)_\Gamma
		=
		\varepsilon(\partial_n\phi_h,r_h)_\Gamma
		+U(J(\phi_h),r_h)_\Gamma,
		& for all $r_h\in Y_h|_\Gamma$, \label{eq:semi-weak-bc-L}\notag\\[1mm]
		\nu(l(\phi_h)(\uhtau-\uw),\bm e_h)_\Gamma
		+\nu(\partial_n\uhtau,\bm e_h)_\Gamma
		-\lambda Q(L(\phi_h)\gradtau\phi_h,\bm e_h)_\Gamma
		=0,
		& for all $\bm e_h\in V_{h,\tau}|_\Gamma$, \label{eq:semi-weak-bc-gnbc}\notag\\[1mm]
		\partial_n\mu_h=0,
		& on $\Gamma$. \label{eq:semi-weak-bc-mu}\notag
	\end{numcases}
\end{subequations}

We now introduce the fully discrete scheme. 
Let $\Delta t>0$ be the time step size and $t^n=n\Delta t$. 
For any sequence $\{a^n\}$, define the BDF2 difference quotient and the second-order extrapolation by
\[
D_t a^{n+1}:=\frac{3a^{n+1}-4a^n+a^{n-1}}{2\Delta t},
\qquad
a^*:=2a^n-a^{n-1}.
\]
For simplicity, we write
\[
H^*=H(\phi_h^*),\qquad
J^*=J(\phi_h^*),
\]
and use the same convention for other explicitly evaluated quantities.

Given
$(\phi_h^n,\mu_h^n,U^n,\uh_h^n,p_h^n,Q^n)$
and the corresponding values at the previous time level, find
\[
(\phi_h^{n+1},\mu_h^{n+1},U^{n+1},
\widetilde{\uh}_h^{n+1},p_h^{n+1},Q^{n+1})
\in
Y_h\times Y_h\times\mathbb R\times V_h\times O_h\times\mathbb R
\]
such that, for all $\omega_h,\eta_h\in Y_h$, $\vh_h\in V_h$, and $q_h\in O_h$,
\begin{subequations}
	\label{eq:scheme}
	\begin{numcases}{}
		(D_t\phi_h^{n+1},\omega_h)
		-
		Q^{n+1}(\uh_h^*\phi_h^*,\nabla\omega_h)
		=
		-M(\nabla\mu_h^{n+1},\nabla\omega_h),
		\label{eq:scheme-phi}\\[1mm]
		\begin{aligned}
			(\mu_h^{n+1},\eta_h)
			=
			\lambda\Big[
			\varepsilon(\nabla\phi_h^{n+1},\nabla\eta_h)
			-\varepsilon(\partial_n\phi_h^{n+1},\eta_h)_\Gamma
			+U^{n+1}(H^*,\eta_h)
			\Big],
		\end{aligned}
		\label{eq:scheme-mu}\\[1mm]
		D_t U^{n+1}
		=
		\frac12(H^*,D_t\phi_h^{n+1})
		+\frac12(J^*,D_t\phi_h^{n+1})_\Gamma,
		\label{eq:scheme-U}\\[1mm]
		\begin{aligned}
			(\frac{3\widetilde{\uh}_h^{n+1}-4{\uh}_h^{n}+{\uh}_h^{n-1}}{2\Delta t},\vh_h)
			&+
			Q^{n+1}((\uh_h^*\cdot\nabla)\uh_h^*,\vh_h)
			+\nu(\nabla\widetilde{\uh}_h^{n+1},\nabla\vh_h)
			-(p_h^n,\nabla\cdot\vh_h)\\
			&-\nu(\partial_n\utildehtaup,\vhtau)_\Gamma
			+Q^{n+1}(\phi_h^*\nabla\mu_h^{*},\vh_h)
			=0,
		\end{aligned}
		\label{eq:scheme-u}\\[1mm]
		\begin{aligned}
			D_t Q^{n+1}
			={}&
			-((\uh_h^*\phi_h^*),\nabla\mu_h^{n+1})
			+(\phi_h^*\nabla\mu_h^*,\widetilde{\uh}_h^{n+1})
			+((\uh_h^*\cdot\nabla)\uh_h^*,\widetilde{\uh}_h^{n+1})
			\\
			&-\lambda
			(\widetilde{L}_h^*\gradtau\phi_h^*,\utildehtaup)_\Gamma
			+\lambda
			(\widetilde{L}_h^{n+1},\tiluhtaus\cdot\gradtau\phi_h^*)_\Gamma.
		\end{aligned}
		\label{eq:scheme-Q}\\[1mm]
			(\nabla\cdot\widetilde{\uh}_h^{n+1},q_h)
		+\frac{2\Delta t}{3}
		(\nabla(p_h^{n+1}-p_h^n),\nabla q_h)
		=0,
		\label{eq:scheme-corr}\\[1mm]
			\uh_h^{n+1}
			=
			\widetilde{\uh}_h^{n+1}
			-\frac{2\Delta t}{3}\nabla(p_h^{n+1}-p_h^n).
			\label{eq:corr-strong}
	\end{numcases}
\end{subequations}
The boundary terms in \eqref{eq:scheme-Q} are chosen to cancel the corresponding coupling terms arising from the discrete GNBC and DCLC in the energy estimate.

The boundary conditions on $\Gamma$ are discretized as follows: for all
$m_h,r_h\in Y_h|\Gamma$ and $\bm e_h\in V_{h,\tau}|\Gamma$,
\begin{subequations}
	\label{eq:scheme-bc}
	\begin{numcases}{}
		(D_t\phi_h^{n+1},m_h)_\Gamma
		+
		Q^{n+1}(\tiluhtaus\cdot\gradtau\phi_h^*,m_h)_\Gamma
		=
		-\gamma(\tilde{L}_h^{n+1},m_h)_\Gamma,
		\label{eq:scheme-bc-dclc}\\[1mm]
		(\tilde{L}_h^{n+1},r_h)_\Gamma
		=
		\varepsilon(\partial_n\phi_h^{n+1},r_h)_\Gamma
		+U^{n+1}(J^*,r_h)_\Gamma,
	  \label{eq:scheme-bc-L}\\[1mm]
		\begin{aligned}
			&\nu(l(\phi_h^n)
			\left(\utildehtaup-\uw\right),\bm e_h)_\Gamma
			+\nu(\partial_n\utildehtaup,\bm e_h)_\Gamma
			-\lambda Q^{n+1}(\tilde{L}_h^{*}\gradtau\phi_h^*,\bm e_h)_\Gamma
			=0,
		\end{aligned}
		 \label{eq:scheme-bc-gnbc}\\[1mm]
		\partial_n\mu_h^{n+1}=0. \label{eq:scheme-bc-mu}
	\end{numcases}
\end{subequations}

\begin{remark}[Linearity and decoupling]
Since the scheme is linear for $U^{n+1}$ and $Q^{n+1}$,
the solution can be expressed in terms of component functions
obtained from linear elliptic subproblems. The scalar coefficients
are then determined from linear algebraic equations.
Thus, each time step requires only solving linear elliptic problems and updating the scalar auxiliary variables, as detailed in Section~\ref{sec:implementation}.
\end{remark}

\begin{remark}
	The pressure-correction step decouples the pressure and velocity
	computations; see \cite{guermond2006} for an overview.
	Combining \eqref{eq:scheme-corr} and \eqref{eq:corr-strong}
	and integrating by parts, we obtain
	\begin{align}
		(\uh_h^{n+1},\nabla q_h)
		&=(\widetilde{\uh}_h^{n+1},\nabla q_h)
		-\frac{2\Delta t}{3}
		(\nabla(p_h^{n+1}-p_h^n),\nabla q_h)\nonumber\\
		&=(\widetilde{\uh}_h^{n+1},\nabla q_h)
		+(\nabla\cdot\widetilde{\uh}_h^{n+1},q_h)=0,
		\qquad \forall q_h\in O_h, 
		\label{no-diver}
	\end{align}
	Here, the boundary integral vanishes by the no-penetration
	condition on the solid walls and the periodic boundary conditions
	on the remaining boundaries.
\end{remark}

\section{Unconditional energy stability}
\label{sec:stability}

\begin{theorem}[Unconditional energy stability]\label{Thm5.1}
	Assume that $\uw=0$. Then the scheme \eqref{eq:scheme}--\eqref{eq:scheme-bc} satisfies the following discrete energy inequality:
	\begin{equation}
		\mathcal E_h^{n+1}-\mathcal E_h^n \le
		-\Delta t\,\nu\|\nabla\widetilde{\uh}_h^{n+1}\|^2
		-\Delta t\,M\|\nabla\mu_h^{n+1}\|^2
		-\Delta t\,\lambda\gamma\|\tilde{L}_h^{n+1}\|_G^2
		-\Delta t\,\nu\|\sqrt{l(\phi_h^n)}\,\utildehtaup\|_\Gamma^2
		\le0,
	\end{equation}
	where the discrete energy $\mathcal E_h^{n+1}$ is defined as
	\begin{equation*}
		\begin{aligned}
			\mathcal E_h^{n+1}:={}&
			\frac14\left(\|\uh_h^{n+1}\|^2+\|2\uh_h^{n+1}-\uh_h^n\|^2\right)
			+\frac{\lambda\varepsilon}{4}\left(\|\nabla\phi_h^{n+1}\|^2+\|2\nabla\phi_h^{n+1}-\nabla\phi_h^n\|^2\right)\\
			&
			+\frac{\lambda}{2}\left(|U^{n+1}|^2+|2U^{n+1}-U^n|^2\right)+\frac14\left(|Q^{n+1}|^2+|2Q^{n+1}-Q^n|^2\right)
			+\frac{\Delta t^2}{3}\|\nabla p_h^{n+1}\|^2-\lambda A.
		\end{aligned}
	\end{equation*}
	Consequently, the scheme is unconditionally energy stable, i.e., no time-step restriction is required for the above energy estimate.
\end{theorem}
\begin{proof} 
	Taking  $\vh_h=2\Delta t\,\widetilde{\uh}_h^{n+1}$ in \eqref{eq:scheme-u}, we obtain
	\begin{align}
		&
		(3\widetilde{\uh}_h^{n+1}-4\uh_h^n+\uh_h^{n-1},\widetilde{\uh}_h^{n+1})
		+2\Delta t\,\nu\|\nabla\widetilde{\uh}_h^{n+1}\|^2
		-2\Delta t\,\nu(\partial_n\utildehtaup,\utildehtaup)_\Gamma
		\notag\\
		&\quad
		-2\Delta t(p_h^n,\nabla\cdot\widetilde{\uh}_h^{n+1})
		+2\Delta t Q^{n+1}
		((\uh_h^*\cdot\nabla)\uh_h^*,\widetilde{\uh}_h^{n+1})\
		+2\Delta t Q^{n+1}
		(\phi_h^*\nabla\mu_h^{*},\widetilde{\uh}_h^{n+1})
		=0.
		\label{eq:proof-u-raw}
	\end{align}
	Taking the $L^2$ inner product of \eqref{eq:corr-strong}  with $\uh_h^{n+1}$, we derive 
	\begin{equation}
		(\widetilde{\uh}_h^{n+1}-{\uh}_h^{n+1},{\uh}_h^{n+1})=\frac{2\Delta t}{3}(\nabla (p_h^{n+1}-p_h^{n}),\uh_h^{n+1})
		=0,
		\label{eq:proof-corr-reformulated}
	\end{equation}
	 and 
		\begin{equation}
	(3{\uh}_h^{n+1}-4\uh_h^n+\uh_h^{n-1},\widetilde{\uh}_h^{n+1}-{\uh}_h^{n+1})=	(3{\uh}_h^{n+1}-4\uh_h^n+\uh_h^{n-1},\frac{2\Delta t}{3}(\nabla (p_h^{n+1}-p_h^{n}))=0,
		\label{eq:proof-corr-reformulated1}
	\end{equation}
	where \eqref{no-diver} is used.
	According to \eqref{eq:proof-corr-reformulated}--\eqref{eq:proof-corr-reformulated1}, we derive the equality
	\begin{align}
		&(3\widetilde{\uh}_h^{n+1}-4\uh_h^n+\uh_h^{n-1},\widetilde{\uh}_h^{n+1})
		\notag\\
		=&(3\widetilde{\uh}_h^{n+1}-3\uh_h^{n+1},\widetilde{\uh}_h^{n+1})
		+(3\uh_h^{n+1}-4\uh_h^n+\uh_h^{n-1},\widetilde{\uh}_h^{n+1})
		\notag\\
		=&(3\widetilde{\uh}_h^{n+1}-3\uh_h^{n+1},\widetilde{\uh}_h^{n+1}+\uh_h^{n+1})
		+(3\uh_h^{n+1}-4\uh_h^{n}+\uh_h^{n-1},\uh_h^{n+1})
		\notag\\
		=&\frac12\Bigl(
		\|\uh_h^{n+1}\|^2-\|\uh_h^n\|^2
		+\|2\uh_h^{n+1}-\uh_h^n\|^2
		-\|2\uh_h^n-\uh_h^{n-1}\|^2
		+\|\uh_h^{n+1}-2\uh_h^n+\uh_h^{n-1}\|^2
		\Bigr)
		\notag\\
		&\quad+3\left(
		\|\widetilde{\uh}_h^{n+1}\|^2
		-\|\uh_h^{n+1}\|^2
		\right).
		\label{eq:proof-velocity-bdf2}
	\end{align}
	Equation~\eqref{eq:corr-strong} can be rewritten as
	\begin{equation}
		\frac{3}{2\Delta t}\uh_h^{n+1}+\nabla p_h^{n+1}
		=
		\frac{3}{2\Delta t}\widetilde{\uh}_h^{n+1}+\nabla p_h^n.
		\label{eq:proof-corr-reformulated2}
	\end{equation}
	Taking the $L^2$ inner product of the above equation with itself and multiplying the result with $2\Delta t^2/3$, we get
	\begin{align}
		\frac32\left(
		\|\uh_h^{n+1}\|^2-\|\widetilde{\uh}_h^{n+1}\|^2
		\right)
		+\frac{2\Delta t^2}{3}\left(
		\|\nabla p_h^{n+1}\|^2-\|\nabla p_h^n\|^2
		\right)
		=2\Delta t(\widetilde{\uh}_h^{n+1},\nabla p_h^n).
		\label{eq:proof-pressure-increment}
	\end{align}
	Rewrite \eqref{eq:corr-strong} as
	\begin{equation*}
		\uh_h^{n+1}-\widetilde{\uh}_h^{n+1}=-\frac{2\Delta t}{3}
		\nabla p_h^{n+1}+\frac{2\Delta t}{3}
		\nabla p_h^{n}.
	\end{equation*}
	By taking the $L^2$ inner product of the above equation with $\frac{3}{2}\uh_h^{n+1}$
	and using \eqref{no-diver}, we obtain
	\begin{equation}
		\frac32\left(
		\|\uh_h^{n+1}\|^2
		-\|\widetilde{\uh}_h^{n+1}\|^2
		+\|\uh_h^{n+1}-\widetilde{\uh}_h^{n+1}\|^2
		\right)=0.
		\label{eq:proof-corr-inner-product}
	\end{equation}
	Combining \eqref{eq:proof-u-raw},
	\eqref{eq:proof-velocity-bdf2},
	\eqref{eq:proof-pressure-increment}, and
	\eqref{eq:proof-corr-inner-product}, we obtain
	\begin{align}
		&
		\frac12
		\left(
		\|\uh_h^{n+1}\|^2-\|\uh_h^n\|^2
		+\|2\uh_h^{n+1}-\uh_h^n\|^2
		-\|2\uh_h^n-\uh_h^{n-1}\|^2
		\right)	+\frac12\|\uh_h^{n+1}-2\uh_h^n+\uh_h^{n-1}\|^2
		\notag\\
		&\quad
		+\frac32\|\uh_h^{n+1}-\widetilde{\uh}_h^{n+1}\|^2
		+\frac{2\Delta t^2}{3}
		\left(
		\|\nabla p_h^{n+1}\|^2
		-\|\nabla p_h^n\|^2
		\right)
		+2\Delta t\,\nu\|\nabla\widetilde{\uh}_h^{n+1}\|^2
		\notag\\
		&\quad
		-2\Delta t\,\nu(\partial_n\utildehtaup,\utildehtaup)_\Gamma
		+2\Delta t Q^{n+1}
		((\uh_h^*\cdot\nabla)\uh_h^*,\widetilde{\uh}_h^{n+1})
		+2\Delta t Q^{n+1}
		(\phi_h^*\nabla\mu_h^{*},\widetilde{\uh}_h^{n+1})
		=0.
		\label{eq:proof-u-est}
	\end{align}
	
	Next, we consider the phase-field part. 
	Taking 
	$w_h=2\Delta t\,\mu_h^{n+1}$ in \eqref{eq:scheme-phi} and using the boundary condition
	$\partial_n\mu_h^{n+1}=0$ on $\Gamma$, we obtain
	\begin{align}
		(3\phi_h^{n+1}-4\phi_h^n+\phi_h^{n-1},\mu_h^{n+1})
		-2\Delta t Q^{n+1}
		(\uh_h^*\phi_h^*,\nabla\mu_h^{n+1})
		+2\Delta t M\|\nabla\mu_h^{n+1}\|^2
		=0.
		\label{eq:proof-phase}
	\end{align}
	Taking
	$\eta=-(3\phi_h^{n+1}-4\phi_h^n+\phi_h^{n-1})$ in \eqref{eq:scheme-mu}, we obtain
	\begin{align}
		&-(\mu_h^{n+1},3\phi_h^{n+1}-4\phi_h^n+\phi_h^{n-1})
		+\lambda\varepsilon
		(\nabla\phi_h^{n+1},\nabla\!\left(
		3\phi_h^{n+1}-4\phi_h^n+\phi_h^{n-1}
		\right))
		\notag\\
		&\quad
		-\lambda\varepsilon
		(\partial_n\phi_h^{n+1},3\phi_h^{n+1}-4\phi_h^n+\phi_h^{n-1})_\Gamma
		+\lambda U^{n+1}
		(H^*,3\phi_h^{n+1}-4\phi_h^n+\phi_h^{n-1})
		=0.
		\label{eq:proof-mu-raw}
	\end{align}
	Multiplying \eqref{eq:scheme-U} by
	$2\lambda U^{n+1}$, we obtain
	\begin{align}
		&
		\lambda
		\left(
		|U^{n+1}|^2-|U^n|^2
		+|2U^{n+1}-U^n|^2
		-|2U^n-U^{n-1}|^2
		+|U^{n+1}-2U^n+U^{n-1}|^2
		\right)
		\notag\\
		&
		=
		\lambda U^{n+1}
		(H^*,3\phi_h^{n+1}-4\phi_h^n+\phi_h^{n-1})
		+\lambda U^{n+1}
		(J^*,3\phi_h^{n+1}-4\phi_h^n+\phi_h^{n-1})_\Gamma.
		\label{eq:proof-U}
	\end{align}
	Multiplying \eqref{eq:scheme-Q} by
	$2\Delta t Q^{n+1}$, we obtain
	\begin{align}
		&
		\frac12
		\left(
		|Q^{n+1}|^2-|Q^n|^2
		+|2Q^{n+1}-Q^n|^2
		-|2Q^n-Q^{n-1}|^2
		+|Q^{n+1}-2Q^n+Q^{n-1}|^2
		\right)
		\notag\\
		&
		=
		-2\Delta t Q^{n+1}
		((\uh_h^*\phi_h^*),\nabla\mu_h^{n+1})
		+2\Delta t Q^{n+1}
		(\phi_h^*\nabla\mu_h^{*},\widetilde{\uh}_h^{n+1})
		+2\Delta t Q^{n+1}
		((\uh_h^*\cdot\nabla)\uh_h^*,\widetilde{\uh}_h^{n+1})
		\notag\\
		&\quad
		-2\Delta t\lambda Q^{n+1}
		(\tilde{L}_h^{*}\gradtau\phi_h^*,\utildehtaup)_\Gamma
		+2\Delta t\lambda Q^{n+1}
		(\tilde{L}_h^{n+1},\tiluhtaus\cdot\gradtau\phi_h^*)_\Gamma.
		\label{eq:proof-Q}
	\end{align}
	For the boundary terms, taking
	$\bm e_h=2\Delta t\,\utildehtaup$ in
	\eqref{eq:scheme-bc-gnbc} and using $\uw=0$, we obtain
	\begin{equation}
		\label{eq:proof-gnbc}
		2\Delta t\nu(\partial_n\utildehtaup,\utildehtaup)_\Gamma
		=
		2\Delta t\lambda Q^{n+1}
		(\tilde{L}_h^{*}\gradtau\phi_h^*,\utildehtaup)_\Gamma
		-2\Delta t\nu\|\sqrt{l(\phi_h^n)}\,\utildehtaup\|_\Gamma^2.
	\end{equation}
	For the phase-field boundary term, 
	\begin{align}
		\lambda\varepsilon
		(\partial_n\phi_h^{n+1},3\phi_h^{n+1}-4\phi_h^n+\phi_h^{n-1})_\Gamma
		=\lambda(\tilde{L}_h^{n+1}-J^*U^{n+1},3\phi_h^{n+1}-4\phi_h^n+\phi_h^{n-1})_\Gamma,
		\label{eq:proof-L}
	\end{align}
	where
	\begin{align}
		\label{eq:proof-dclc}
		&\lambda(\tilde{L}_h^{n+1},3\phi_h^{n+1}-4\phi_h^n+\phi_h^{n-1})_\Gamma\notag\\
		=&2\Delta t\lambda(\tilde{L}_h^{n+1},-\gamma \tilde{L}_h^{n+1}
		-Q^{n+1}\tiluhtaus\cdot\gradtau\phi_h^*)_\Gamma
		\notag\\
		=&-2\Delta t\lambda\gamma\|\tilde{L}_h^{n+1}\|_\Gamma^2
		-2\Delta t\lambda Q^{n+1}
		(\tilde{L}_h^{n+1},\tiluhtaus\cdot\gradtau\phi_h^*)_\Gamma.
	\end{align}
	
	Summing \eqref{eq:proof-u-est}--
	\eqref{eq:proof-Q}, and using the boundary identities
	\eqref{eq:proof-gnbc}--\eqref{eq:proof-dclc}, we obtain
	\begin{align}
		&
		\frac12
		\left(
		\|\uh_h^{n+1}\|^2-\|\uh_h^n\|^2
		+\|2\uh_h^{n+1}-\uh_h^n\|^2
		-\|2\uh_h^n-\uh_h^{n-1}\|^2
		\right)+\frac{2\Delta t^2}{3}
		\left(
		\|\nabla p_h^{n+1}\|^2-\|\nabla p_h^n\|^2
		\right)
		\notag\\
		&\quad
		+\frac{\lambda\varepsilon}{2}
		\left(
		\|\nabla\phi_h^{n+1}\|^2-\|\nabla\phi_h^n\|^2
		+\|2\nabla\phi_h^{n+1}-\nabla\phi_h^n\|^2
		-\|2\nabla\phi_h^n-\nabla\phi_h^{n-1}\|^2
		\right)
		\notag\\
		&\quad
		+\lambda
		\left(
		|U^{n+1}|^2-|U^n|^2
		+|2U^{n+1}-U^n|^2
		-|2U^n-U^{n-1}|^2
		\right)
		\notag\\
		&\quad
		+\frac12
		\left(
		|Q^{n+1}|^2-|Q^n|^2
		+|2Q^{n+1}-Q^n|^2
		-|2Q^n-Q^{n-1}|^2
		\right)
		\notag\\
		&\quad
		+\frac12\|\uh_h^{n+1}-2\uh_h^n+\uh_h^{n-1}\|^2
		+\frac32\|\uh_h^{n+1}-\widetilde{\uh}_h^{n+1}\|^2
		+\frac{\lambda\varepsilon}{2}
		\|\nabla(\phi_h^{n+1}-2\phi_h^n+\phi_h^{n-1})\|^2
		\notag\\
		&\quad
		+\lambda|U^{n+1}-2U^n+U^{n-1}|^2
		+\frac12|Q^{n+1}-2Q^n+Q^{n-1}|^2
		\notag\\
		&\le
		-2\Delta t\nu\|\nabla\widetilde{\uh}_h^{n+1}\|^2
		-2\Delta t M\|\nabla\mu_h^{n+1}\|^2
		-2\Delta t\nu\|\sqrt{l(\phi_h^n)}\,\utildehtaup\|_\Gamma^2
		-2\Delta t\lambda\gamma\|\tilde{L}_h^{n+1}\|_\Gamma^2.
		\label{eq:proof-final-raw}
	\end{align}
	Then we obtain the desired result.
\end{proof}
\begin{remark}
	The auxiliary-variable contributions to $\mathcal E_h^{n+1}$
	are nonnegative quadratic terms. Consequently,
	\[
	\mathcal E_h^{n+1}\ge -\lambda A.
	\]
	This differs from the step-by-step SAV energy in
	\cite{gao2025}, which contains linear combinations
	of the bulk and wall energy variables whose
	nonnegativity is not guaranteed; see Remark~3.5 therein.
\end{remark}
\section{Decoupled implementation}
\label{sec:implementation}

In this section, we present the solution procedure for the fully discrete scheme, in which the subproblems are solved sequentially.

\subsection{Splitting with respect to the ZEC variable}

We first use the scalar ZEC variable $Q^{n+1}$ to split the phase-field variables and the unified SAV variable into the form
\begin{subequations}
	\label{eq:splitting-Q}
	\begin{numcases}{}
		\phi_h^{n+1}
		=
		\phi_{1h}^{n+1}
		+
		Q^{n+1}\phi_{2h}^{n+1},
		\label{eq:splitting-Q-phi}\\
		\mu_h^{n+1}
		=
		\mu_{1h}^{n+1}
		+
		Q^{n+1}\mu_{2h}^{n+1},
		\label{eq:splitting-Q-mu}\\
		U^{n+1}
		=
		U_1^{n+1}
		+
		Q^{n+1}U_2^{n+1}.
		\label{eq:splitting-Q-U}
	\end{numcases}
\end{subequations}
Substituting \eqref{eq:splitting-Q} into the discrete phase equation and chemical-potential equation, then collecting the terms independent of $Q^{n+1}$ and the terms proportional to $Q^{n+1}$, we obtain two subproblems.

The first subproblem is: find $(\phi_{1h}^{n+1},\mu_{1h}^{n+1})\in Y_h\times Y_h$ such that for all $\omega_h,\eta_h\in Y_h$,
\begin{subequations}
	\label{eq:subproblem-phi1}
	\begin{numcases}{}
		\dfrac{3}{2M\Delta t}(\phi_{1h}^{n+1},\omega_h)
		+
		(\nabla\mu_{1h}^{n+1},\nabla\omega_h)
		=
		\dfrac{1}{2M\Delta t}
		(4\phi_h^n-\phi_h^{n-1},\omega_h),
		\label{eq:subproblem-phi1-a}\\[2mm]
		\begin{aligned}
			(\mu_{1h}^{n+1},\eta_h)
			=
			\lambda\Big[
			&\varepsilon(\nabla\phi_{1h}^{n+1},\nabla\eta_h)
			-\varepsilon(\partial_n\phi_{1h}^{n+1},\eta_h)_\Gamma
			+U_1^{n+1}(H^*,\eta_h)
			\Big].
		\end{aligned}
		\label{eq:subproblem-phi1-b}
	\end{numcases}
\end{subequations}

The second subproblem is: find $(\phi_{2h}^{n+1},\mu_{2h}^{n+1})\in Y_h\times Y_h$ such that
\begin{subequations}
	\label{eq:subproblem-phi2}
	\begin{numcases}{}
		\dfrac{3}{2M\Delta t}(\phi_{2h}^{n+1},\omega_h)
		+
		(\nabla\mu_{2h}^{n+1},\nabla\omega_h)
		=
		\dfrac{1}{M}(\uh_h^*\phi_h^*,\nabla\omega_h),
		\label{eq:subproblem-phi2-a}\\[2mm]
		\begin{aligned}
			(\mu_{2h}^{n+1},\eta_h)
			=
			\lambda\Big[
			&\varepsilon(\nabla\phi_{2h}^{n+1},\nabla\eta_h)
			-\varepsilon(\partial_n\phi_{2h}^{n+1},\eta_h)_\Gamma
			+U_2^{n+1}(H^*,\eta_h)
			\Big].
		\end{aligned}
		\label{eq:subproblem-phi2-b}
	\end{numcases}
\end{subequations}
The corresponding boundary equations are obtained from \eqref{eq:scheme-bc-dclc}--\eqref{eq:scheme-bc-L}. 
For the first component, we have
\begin{equation}
	\label{eq:bc-phi1}
	\frac{3}{2\Delta t}(\phi_{1h}^{n+1},m_h)_\Gamma
	+
	\gamma\varepsilon(\partial_n\phi_{1h}^{n+1},m_h)_\Gamma
	+
	\gamma U_1^{n+1}(J^*,m_h)_\Gamma
	=
	\frac{1}{2\Delta t}(4\phi_h^n-\phi_h^{n-1},m_h)_\Gamma,
\end{equation}
for all $m_h\in Y_h|_\Gamma$. 
For the second component, we have
\begin{equation}
	\label{eq:bc-phi2}
	\frac{3}{2\Delta t}(\phi_{2h}^{n+1},m_h)_\Gamma
	+
	\gamma\varepsilon(\partial_n\phi_{2h}^{n+1},m_h)_\Gamma
	+
	\gamma U_2^{n+1}(J^*,m_h)_\Gamma
	=
	-(\tiluhtaus\cdot\gradtau\phi_h^*,m_h)_\Gamma.
\end{equation}
At this stage, the two subproblems still contain the scalar coefficients $U_1^{n+1}$ and $U_2^{n+1}$. 
Then we remove this remaining scalar coupling.

\subsection{Splitting with respect to the unified SAV variable}
Since $U_1^{n+1}$ and $U_2^{n+1}$ are scalar unknowns, each of the two subproblems above can be further split. 
For the first subproblem, we write
\begin{subequations}
	\label{eq:splitting-U1}
	\begin{numcases}{}
		\phi_{1h}^{n+1}
		=
		\phi_{11h}^{n+1}
		+
		U_1^{n+1}\phi_{12h}^{n+1},
		\label{eq:splitting-U1-phi}\\
		\mu_{1h}^{n+1}
		=
		\mu_{11h}^{n+1}
		+
		U_1^{n+1}\mu_{12h}^{n+1}.
		\label{eq:splitting-U1-mu}
	\end{numcases}
\end{subequations}
The component independent of $U_1^{n+1}$ is obtained by finding $(\phi_{11h}^{n+1},\mu_{11h}^{n+1})\in Y_h\times Y_h$ such that
\begin{subequations}
	\label{eq:phi11-system}
	\begin{numcases}{}
		\dfrac{3}{2M\Delta t}(\phi_{11h}^{n+1},\omega_h)
		+
		(\nabla\mu_{11h}^{n+1},\nabla\omega_h)
		=
		\dfrac{1}{2M\Delta t}
		(4\phi_h^n-\phi_h^{n-1},\omega_h),
		\label{eq:phi11-a}\\[2mm]
		(\mu_{11h}^{n+1},\eta_h)
		=
		\lambda\varepsilon(\nabla\phi_{11h}^{n+1},\nabla\eta_h)
		-\lambda\varepsilon(\partial_n\phi_{11h}^{n+1},\eta_h)_\Gamma,
		\label{eq:phi11-b}
	\end{numcases}
\end{subequations}
with boundary condition
\begin{equation}
	\label{eq:phi11-bc}
	\frac{3}{2\Delta t}(\phi_{11h}^{n+1},m_h)_\Gamma
	+
	\gamma\varepsilon(\partial_n\phi_{11h}^{n+1},m_h)_\Gamma
	=
	\frac{1}{2\Delta t}(4\phi_h^n-\phi_h^{n-1},m_h)_\Gamma.
\end{equation}
The component proportional to $U_1^{n+1}$ is obtained by finding $(\phi_{12h}^{n+1},\mu_{12h}^{n+1})\in Y_h\times Y_h$ such that
\begin{subequations}
	\label{eq:phi12-system}
	\begin{numcases}{}
		\dfrac{3}{2M\Delta t}(\phi_{12h}^{n+1},\omega_h)
		+
		(\nabla\mu_{12h}^{n+1},\nabla\omega_h)
		=0,
		\label{eq:phi12-a}\\[2mm]
		(\mu_{12h}^{n+1},\eta_h)
		=
		\lambda\varepsilon(\nabla\phi_{12h}^{n+1},\nabla\eta_h)
		-\lambda\varepsilon(\partial_n\phi_{12h}^{n+1},\eta_h)_\Gamma
		+\lambda(H^*,\eta_h),
		\label{eq:phi12-b}
	\end{numcases}
\end{subequations}
with boundary condition
\begin{equation}
	\label{eq:phi12-bc}
	\frac{3}{2\Delta t}(\phi_{12h}^{n+1},m_h)_\Gamma
	+
	\gamma\varepsilon(\partial_n\phi_{12h}^{n+1},m_h)_\Gamma
	+
	\gamma(J^*,m_h)_\Gamma
	=0.
\end{equation}

The second subproblem is treated in the same way. 
We write
\begin{subequations}
	\label{eq:splitting-U2}
	\begin{numcases}{}
		\phi_{2h}^{n+1}
		=
		\phi_{21h}^{n+1}
		+
		U_2^{n+1}\phi_{22h}^{n+1},
		\label{eq:splitting-U2-phi}\\
		\mu_{2h}^{n+1}
		=
		\mu_{21h}^{n+1}
		+
		U_2^{n+1}\mu_{22h}^{n+1}.
		\label{eq:splitting-U2-mu}
	\end{numcases}
\end{subequations}
The two pairs $(\phi_{21h}^{n+1},\mu_{21h}^{n+1})$ and $(\phi_{22h}^{n+1},\mu_{22h}^{n+1})$ are computed from the same left-hand operator as in \eqref{eq:phi11-system}--\eqref{eq:phi12-system}, but with right-hand sides corresponding to \eqref{eq:subproblem-phi2} and \eqref{eq:bc-phi2}. 
Specifically, $(\phi_{21h}^{n+1},\mu_{21h}^{n+1})$ solves
\begin{subequations}
	\label{eq:phi21-system}
	\begin{numcases}{}
		\dfrac{3}{2M\Delta t}(\phi_{21h}^{n+1},\omega_h)
		+
		(\nabla\mu_{21h}^{n+1},\nabla\omega_h)
		=
		\dfrac{1}{M}(\uh_h^*\phi_h^*,\nabla\omega_h),
		\label{eq:phi21-a}\\[2mm]
		(\mu_{21h}^{n+1},\eta_h)
		=
		\lambda\varepsilon(\nabla\phi_{21h}^{n+1},\nabla\eta_h)
		-\lambda\varepsilon(\partial_n\phi_{21h}^{n+1},\eta_h)_\Gamma,
		\label{eq:phi21-b}
	\end{numcases}
\end{subequations}
with
\begin{equation*}
	\label{eq:phi21-bc}
	\frac{3}{2\Delta t}(\phi_{21h}^{n+1},m_h)_\Gamma
	+
	\gamma\varepsilon(\partial_n\phi_{21h}^{n+1},m_h)_\Gamma
	=
	-(\tiluhtaus\cdot\gradtau\phi_h^*,m_h)_\Gamma.
\end{equation*}
The pair $(\phi_{22h}^{n+1},\mu_{22h}^{n+1})$ solves
\begin{subequations}
	\label{eq:phi22-system}
	\begin{numcases}{}
		\dfrac{3}{2M\Delta t}(\phi_{22h}^{n+1},\omega_h)
		+
		(\nabla\mu_{22h}^{n+1},\nabla\omega_h)
		=0,
		\label{eq:phi22-a}\\[2mm]
		(\mu_{22h}^{n+1},\eta_h)
		=
		\lambda\varepsilon(\nabla\phi_{22h}^{n+1},\nabla\eta_h)
		-\lambda\varepsilon(\partial_n\phi_{22h}^{n+1},\eta_h)_\Gamma
		+\lambda(H^*,\eta_h),
		\label{eq:phi22-b}
	\end{numcases}
\end{subequations}
with
\begin{equation*}
	\label{eq:phi22-bc}
	\frac{3}{2\Delta t}(\phi_{22h}^{n+1},m_h)_\Gamma
	+
	\gamma\varepsilon(\partial_n\phi_{22h}^{n+1},m_h)_\Gamma
	+
	\gamma(J^*,m_h)_\Gamma
	=0.
\end{equation*}
Thus the phase-field update is reduced to four linear elliptic subproblems with the same left-hand side structure.

\subsection{Determination of the SAV components}

We now determine $U_1^{n+1}$ and $U_2^{n+1}$. 
The BDF2 discretization of the SAV equation can be written as
\begin{equation}
	\label{eq:U-rewritten}
	U^{n+1}
	=
	\frac12(H^*,\phi_h^{n+1})
	+
	\frac12(J^*,\phi_h^{n+1})_\Gamma
	+
	g^n,
\end{equation}
where the known quantity $g^n$ is given by
\begin{equation*}
	\label{eq:g-n}
	g^n
	=
	\frac{1}{3}(4U^n-U^{n-1})
	-\frac{1}{6}(H^*,4\phi_h^n-\phi_h^{n-1})
	-\frac{1}{6}(J^*,4\phi_h^n-\phi_h^{n-1})_\Gamma.
\end{equation*}
Substituting \eqref{eq:splitting-Q}, \eqref{eq:splitting-U1}, and \eqref{eq:splitting-U2} into \eqref{eq:U-rewritten}, and comparing the coefficients of $1$ and $Q^{n+1}$, we obtain
\begin{subequations}
	\begin{numcases}{}
		U_1^{n+1}
		=
		\dfrac12(H^*,\phi_{11h}^{n+1}+U_1^{n+1}\phi_{12h}^{n+1})
		+
		\dfrac12(J^*,\phi_{11h}^{n+1}+U_1^{n+1}\phi_{12h}^{n+1})_\Gamma
		+
		g^n,
		\label{eq:U1-equation}\notag\\[2mm]
		U_2^{n+1}
		=
		\dfrac12(H^*,\phi_{21h}^{n+1}+U_2^{n+1}\phi_{22h}^{n+1})
		+
		\dfrac12(J^*,\phi_{21h}^{n+1}+U_2^{n+1}\phi_{22h}^{n+1})_\Gamma.\notag
		\label{eq:U2-equation}
	\end{numcases}
\end{subequations}
Hence
\begin{subequations}
	\label{eq:U12-formula}
	\begin{numcases}{}
		\displaystyle
		U_1^{n+1}
		=
		\frac{
			\dfrac12(H^*,\phi_{11h}^{n+1})
			+
			\dfrac12(J^*,\phi_{11h}^{n+1})_\Gamma
			+
			g^n
		}{
			1-\dfrac12(H^*,\phi_{12h}^{n+1})
			-\dfrac12(J^*,\phi_{12h}^{n+1})_\Gamma
		},
		\label{eq:U1-formula}\\[5mm]
		\displaystyle
		U_2^{n+1}
		=
		\frac{
			\dfrac12(H^*,\phi_{21h}^{n+1})
			+
			\dfrac12(J^*,\phi_{21h}^{n+1})_\Gamma
		}{
			1-\dfrac12(H^*,\phi_{22h}^{n+1})
			-\dfrac12(J^*,\phi_{22h}^{n+1})_\Gamma
		}.
		\label{eq:U2-formula}
	\end{numcases}
\end{subequations}
We next verify that the scalar formulas for $U_1^{n+1}$ and
$U_2^{n+1}$ are well defined. Taking $\omega_h=-\mu_{12h}^{n+1}$
in \eqref{eq:phi12-a} and
$\eta_h=\dfrac{3}{2M\Delta t\lambda}\phi_{12h}^{n+1}$ in
\eqref{eq:phi12-b}, we obtain
\begin{align}
	-\frac12(H^*,\phi_{12h}^{n+1})
	=
	\frac{M\Delta t}{3\lambda}\|\nabla\mu_{12h}^{n+1}\|^2
	-\frac{\varepsilon}{2}
	(\partial_n\phi_{12h}^{n+1},\phi_{12h}^{n+1})_\Gamma
	+\frac{\varepsilon}{2}\|\nabla\phi_{12h}^{n+1}\|^2.
	\label{eq:U1-denominator-bulk}
\end{align}
Taking $m_h=\phi_{12h}^{n+1}$ in the boundary equation
\eqref{eq:phi12-bc} and multiplying the resulting identity by
$1/(2\gamma)$, we obtain
\begin{equation}
	-\frac12(J^*,\phi_{12h}^{n+1})_\Gamma
	=
	\frac{3}{4\gamma\Delta t}\|\phi_{12h}^{n+1}\|_\Gamma^2
	+\frac{\varepsilon}{2}
	(\partial_n\phi_{12h}^{n+1},\phi_{12h}^{n+1})_\Gamma.
	\label{eq:U1-denominator-boundary}
\end{equation}
Adding \eqref{eq:U1-denominator-bulk} and
\eqref{eq:U1-denominator-boundary}, we obtain
\begin{align}
	&1-\frac12(H^*,\phi_{12h}^{n+1})
	-\frac12(J^*,\phi_{12h}^{n+1})_\Gamma
	\notag\\
	=&
	1+\frac{M\Delta t}{3\lambda}\|\nabla\mu_{12h}^{n+1}\|^2
	+\frac{\varepsilon}{2}\|\nabla\phi_{12h}^{n+1}\|^2
	+\frac{3}{4\gamma\Delta t}\|\phi_{12h}^{n+1}\|_\Gamma^2
	>0.
	\label{eq:U-denominator-positive}\notag
\end{align}
The pair $(\phi_{22h}^{n+1},\mu_{22h}^{n+1})$ satisfies the same
bulk and boundary equations. Repeating the above
argument with $(\phi_{22h}^{n+1},\mu_{22h}^{n+1})$ therefore proves
that the denominator in \eqref{eq:U2-formula} is also strictly
positive. Hence, $U_1^{n+1}$ and $U_2^{n+1}$ are uniquely determined.

After $U_1^{n+1}$ and $U_2^{n+1}$ are known, we recover
\begin{subequations}
	\label{eq:recover-phi12}
	\begin{numcases}{}
		\phi_{1h}^{n+1}
		=
		\phi_{11h}^{n+1}
		+
		U_1^{n+1}\phi_{12h}^{n+1},
		\\
		\mu_{1h}^{n+1}
		=
		\mu_{11h}^{n+1}
		+
		U_1^{n+1}\mu_{12h}^{n+1},
		\\
		\phi_{2h}^{n+1}
		=
		\phi_{21h}^{n+1}
		+
		U_2^{n+1}\phi_{22h}^{n+1},
		\\
		\mu_{2h}^{n+1}
		=
		\mu_{21h}^{n+1}
		+
		U_2^{n+1}\mu_{22h}^{n+1}.
	\end{numcases}
\end{subequations}

\subsection{Velocity splitting}

The intermediate velocity is also split with respect to $Q^{n+1}$:
\begin{equation}
	\label{eq:velocity-splitting}
	\widetilde{\uh}_h^{n+1}
	=
	\widetilde{\uh}_{1h}^{n+1}
	+
	Q^{n+1}\widetilde{\uh}_{2h}^{n+1}.
\end{equation}
Substituting \eqref{eq:velocity-splitting} into \eqref{eq:scheme-u} gives two linear elliptic problems for the intermediate-velocity components.

The first velocity component is determined by finding $\widetilde{\uh}_{1h}^{n+1}\in V_h$ such that for all $\vh_h\in V_h$,
\begin{equation}
	\label{eq:u1-system}
	\frac{3}{2\Delta t}(\widetilde{\uh}_{1h}^{n+1},\vh_h)
	+\nu(\nabla\widetilde{\uh}_{1h}^{n+1},\nabla\vh_h)
	-\nu(\partial_n\widetilde{\uh}_{1h,\tau}^{n+1},\vhtau)_\Gamma
	=
	(\frac{4\uh_h^n-\uh_h^{n-1}}{2\Delta t},\vh_h)
	+(p_h^n,\nabla\cdot\vh_h).
\end{equation}
Its boundary condition is
\begin{equation*}
	\label{eq:u1-bc}
	\nu(l(\phi_h^n)\widetilde{\uh}_{1h,\tau}^{n+1},\bm e_h)_\Gamma
	+
	\nu(\partial_n\widetilde{\uh}_{1h,\tau}^{n+1},\bm e_h)_\Gamma
	=0,
	\qquad
	\forall \bm e_h\in V_{h,\tau}|_\Gamma.
\end{equation*}
The second velocity component $\widetilde{\bm u}_{2h}^{n+1}\in V_h$
satisfies the following equation for all $\bm v_h\in V_h$:
\begin{equation}
	\label{eq:u2-system}
	\frac{3}{2\Delta t}(\widetilde{\uh}_{2h}^{n+1},\vh_h)
	+\nu(\nabla\widetilde{\uh}_{2h}^{n+1},\nabla\vh_h)
	-\nu(\partial_n\widetilde{\uh}_{2h,\tau}^{n+1},\vhtau)_\Gamma
	=
	-((\uh_h^*\cdot\nabla)\uh_h^*,\vh_h)
	-(\phi_h^*\nabla\mu_h^{*},\vh_h).
\end{equation}
The corresponding boundary condition is
\begin{equation*}
	\label{eq:u2-bc}
	\nu(l(\phi_h^n)\widetilde{\uh}_{2h,\tau}^{n+1},\bm e_h)_\Gamma
	+
	\nu(\partial_n\widetilde{\uh}_{2h,\tau}^{n+1},\bm e_h)_\Gamma
	=
	\lambda(\tilde{L}_h^{*
	}\gradtau\phi_h^*,\bm e_h),
	\qquad
	\forall \bm e_h\in V_{h,\tau}|_\Gamma.
\end{equation*}
Here
\[
\widetilde{\uh}_{ih,\tau}^{n+1}
=
\widetilde{\uh}_{ih}^{n+1}
-
(\widetilde{\uh}_{ih}^{n+1}\cdot\nh)\nh,
\qquad i=1,2.
\]
Since the left-hand side operators in \eqref{eq:u1-system} and \eqref{eq:u2-system} have the same structure, the same matrix can be reused when the coefficients are treated explicitly.

\subsection{Determination of the ZEC variable}

We next determine $Q^{n+1}$. 
Substituting the decompositions
\[
\mu_h^{n+1}
=
\mu_{1h}^{n+1}
+
Q^{n+1}\mu_{2h}^{n+1},
\qquad
\widetilde{\uh}_h^{n+1}
=
\widetilde{\uh}_{1h}^{n+1}
+
Q^{n+1}\widetilde{\uh}_{2h}^{n+1}
\]
into the ZEC equation \eqref{eq:scheme-Q}, we obtain a scalar equation of the form
\begin{equation}
	\label{eq:Q-scalar}
	\left(
	\frac{3}{2\Delta t}
	-
	\mathcal V_2
	\right)Q^{n+1}
	=
	\frac{4Q^n-Q^{n-1}}{2\Delta t}
	+
	\mathcal V_1,
\end{equation}
where, for $i=1,2$,
\begin{align}
	\mathcal V_i
	={}&
	-(\uh_h^*\phi_h^*,\nabla\mu_{ih}^{n+1})
	+
	(\phi_h^*\nabla\mu_h^{*},\widetilde{\uh}_{ih}^{n+1})
	+
	((\uh_h^*\cdot\nabla)\uh_h^*,\widetilde{\uh}_{ih}^{n+1})
	\notag\\
	&-\lambda(\tilde{L}_h^{*}\gradtau\phi_h^*,\widetilde{\uh}_{ih,\tau}^{n+1})_\Gamma
	+\lambda(\tilde{L}_{ih}^{n+1},\tiluhtaus\cdot\gradtau\phi_h^*)_\Gamma.
	\label{eq:Vi-definition}
\end{align}
We subsequently verify that \eqref{eq:Q-scalar} is solvable by proving that
$\dfrac{3}{2\Delta t}-\mathcal V_2>0$. Taking $\vh_h=
\widetilde{\uh}_{2h}^{n+1}$ in \eqref{eq:u2-system}, we obtain
\begin{align}
	-(\phi_h^*\nabla\mu_h^*,\widetilde{\uh}_{2h}^{n+1})
	-((\uh_h^*\cdot\nabla)\uh_h^*,\widetilde{\uh}_{2h}^{n+1})
	=
	\frac{3}{2\Delta t}\|\widetilde{\uh}_{2h}^{n+1}\|^2
	+\nu\|\nabla\widetilde{\uh}_{2h}^{n+1}\|^2
	-\nu(\partial_n\widetilde{\uh}_{2h,\tau}^{n+1},\widetilde{\uh}_{2h,\tau}^{n+1})_\Gamma.
	\notag
\end{align}
Taking $\omega_h=M\mu_{2h}^{n+1}$ in
\eqref{eq:subproblem-phi2-a} and
$\eta_h=-\dfrac{3}{2\Delta t}\phi_{2h}^{n+1}$ in
\eqref{eq:subproblem-phi2-b}, and combining the resulting two
equalities, we obtain
\begin{align}
	(\uh_h^*\phi_h^*,\nabla\mu_{2h}^{n+1})
	=
	M\|\nabla\mu_{2h}^{n+1}\|^2
	+\frac{3\lambda\varepsilon}{2\Delta t}
	\|\nabla\phi_{2h}^{n+1}\|^2
	-\frac{3\lambda\varepsilon}{2\Delta t}
	(\partial_n\phi_{2h}^{n+1},\phi_{2h}^{n+1})_\Gamma
	+\frac{3\lambda}{2\Delta t}U_2^{n+1}
	(H^*,\phi_{2h}^{n+1}).
	\notag
\end{align}
For the velocity boundary term, we derive
\begin{align}
	-\nu(\partial_n\widetilde{\uh}_{2h}^{n+1},\widetilde{\uh}_{2h}^{n+1})_\Gamma
	&=
	(
	\nu l(\phi_h^n)\widetilde{\uh}_{2h,s}^{n+1}
	-\lambda
	\bigl(
	\varepsilon\partial_n\phi_h^*
	+J^*U^*
	\bigr)\gradtau\phi_h^*,
	{\widetilde{\uh}_{2h}^{n+1}})_\Gamma
	\notag\\
	&=
	\nu(
	l(\phi_h^n)\widetilde{\uh}_{2h,s}^{n+1},
	{\widetilde{\uh}_{2h,\tau}^{n+1}})_\Gamma
	-\lambda(
	(
	\varepsilon\partial_n\phi_h^*
	+J^*U^*
	)\gradtau\phi_h^*,
	{\widetilde{\uh}_{2h,\tau}^{n+1}})_\Gamma.\notag
\end{align}
where the last term cancels the fourth term in
$\mathcal V_2$.
\\
For the phase boundary term, using \eqref{eq:bc-phi2}, we have
\begin{align}
	&-\frac{3\lambda\varepsilon}{2\Delta t}
	(\partial_n\phi_{2h}^{n+1},\phi_{2h}^{n+1})_\Gamma
	\notag\\
	=&
	-\frac{3\lambda}{2\Delta t}
	(
	-\frac{1}{\gamma}\tiluhtaus\cdot\gradtau\phi_h^*
	-\frac{3}{2\gamma\Delta t}\phi_{2h}^{n+1}
	-J^*U_2^{n+1},
	{\phi_{2h}^{n+1}})_\Gamma
	\notag\\
	=&
	\frac{3\lambda}{2\gamma\Delta t}
	(\tiluhtaus\cdot\gradtau\phi_h^*,\phi_{2h}^{n+1})_\Gamma
	+\frac{9\lambda}{4\gamma\Delta t^2}
	\|\phi_{2h}^{n+1}\|_\Gamma^2
	+\frac{3\lambda}{2\Delta t}U_2^{n+1}
	(J^*,\phi_{2h}^{n+1})_\Gamma.
	\label{eq:Q-solvability-phase-boundary}
\end{align}
The last boundary term in $\mathcal V_2$ satisfies
\begin{align}
	&-\lambda\int_{\Gamma}
	\bigl(
	\varepsilon\partial_n\phi_{2h}^{n+1}
	+J^*U_2^{n+1}
	\bigr)
	\tiluhtaus\cdot\gradtau\phi_h^*\,ds
	\notag\\
	&\quad=
	\lambda(
	\varepsilon\partial_n\phi_{2h}^{n+1}
	+J^*U_2^{n+1},\frac{3}{2\Delta t}\phi_{2h}^{n+1}
	+\gamma\varepsilon\partial_n\phi_{2h}^{n+1}
	+\gamma J^*U_2^{n+1})_\Gamma
	\notag\\
	&\quad=
	\frac{3\lambda}{2\Delta t}
	(
	\varepsilon\partial_n\phi_{2h}^{n+1}
	+J^*U_2^{n+1},\phi_{2h}^{n+1})_\Gamma
	+
	\lambda\gamma
	\|
	\varepsilon\partial_n\phi_{2h}^{n+1}
	+J^*U_2^{n+1}\|_\Gamma^2.\notag
\end{align}
Combining this equality with the first term in
\eqref{eq:Q-solvability-phase-boundary}, we obtain
\begin{equation}
	\frac{3\lambda}{2\Delta t}
	(
	\gamma\varepsilon\partial_n\phi_{2h}^{n+1}
	+\gamma J^*U_2^{n+1}
	+\tiluhtaus\cdot\gradtau\phi_h^*,\phi_{2h}^{n+1})_\Gamma
	=
	-\frac{9\lambda}{4\Delta t^2\gamma}
	\|\phi_{2h}^{n+1}\|_\Gamma^2,
	\label{eq:Q-solvability-boundary-cancellation}\notag
\end{equation}
which cancels the second term in
\eqref{eq:Q-solvability-phase-boundary}.
Furthermore, it follows from \eqref{eq:U2-equation} that
\begin{align*}
	&\frac{3\lambda}{2\Delta t}U_2^{n+1}
	\int_{\Omega}H^*\phi_{2h}^{n+1}\,dx
	+
	\frac{3\lambda}{2\Delta t}U_2^{n+1}
	\int_{\Gamma}J^*\phi_{2h}^{n+1}\,ds
	=
	\frac{3\lambda}{\Delta t}
{|U_2^{n+1}|}^2
	\geq 0.
\end{align*}
Therefore, we obtain $-\mathcal V_2\geq 0$,
which implies that \eqref{eq:Q-scalar} is solvable.

\subsection{Recovery of the full solution and pressure correction}
Once $Q^{n+1}$ has been computed, the phase-field variables and the unified SAV variable can be recovered from
\begin{subequations}
	\label{eq:recover-final-phase}
	\begin{numcases}{}
		\phi_h^{n+1}
		=
		\phi_{1h}^{n+1}
		+
		Q^{n+1}\phi_{2h}^{n+1},
		\\
		\mu_h^{n+1}
		=
		\mu_{1h}^{n+1}
		+
		Q^{n+1}\mu_{2h}^{n+1},
		\\
		U^{n+1}
		=
		U_1^{n+1}
		+
		Q^{n+1}U_2^{n+1}.
	\end{numcases}
\end{subequations}
Similarly, the intermediate velocity is recovered from
\begin{equation}
	\label{eq:recover-u-tilde}
	\widetilde{\uh}_h^{n+1}
	=
	\widetilde{\uh}_{1h}^{n+1}
	+
	Q^{n+1}\widetilde{\uh}_{2h}^{n+1},
\end{equation}
Finally, the pressure $p_h^{n+1}$ and corrected velocity
$\uh_h^{n+1}$ are obtained from
\eqref{eq:scheme-corr}--\eqref{eq:corr-strong}.

The whole implementation is summarized below.

\begin{enumerate}[label=\textbf{Step \arabic*.}]
	\item Compute the extrapolated variables $\phi_h^*$, $\uh_h^*$, evaluate $H^*$, $J^*$, $l(\phi_h^n)$, and other explicitly treated coefficients.
	
	\item Solve the four phase-field component systems \eqref{eq:phi11-system}, \eqref{eq:phi12-system}, \eqref{eq:phi21-system}, and \eqref{eq:phi22-system}, together with their boundary equations.
	
	\item Compute $U_1^{n+1}$ and $U_2^{n+1}$ from \eqref{eq:U1-formula}--\eqref{eq:U2-formula}. 
	Then recover $(\phi_{1h}^{n+1},\mu_{1h}^{n+1})$ and $(\phi_{2h}^{n+1},\mu_{2h}^{n+1})$ from \eqref{eq:recover-phi12}.
	
	\item Solve the two intermediate-velocity component systems \eqref{eq:u1-system} and \eqref{eq:u2-system}, together with their GNBC boundary equations.
	
	\item Determine $Q^{n+1}$ from the scalar equation \eqref{eq:Q-scalar}.
	
	\item Assemble $\phi_h^{n+1}$, $\mu_h^{n+1}$, $U^{n+1}$, and $\widetilde{\uh}_h^{n+1}$ using \eqref{eq:recover-final-phase} and \eqref{eq:recover-u-tilde}.
	
	\item Perform the pressure-correction step \eqref{eq:scheme-corr}--\eqref{eq:corr-strong} to obtain $\uh_h^{n+1}$ and $p_h^{n+1}$.
\end{enumerate}

This procedure shows that one only requires solving linear elliptic subproblems and scalar algebraic equations at each time step. 
No nonlinear iteration is needed. 
Moreover, for fixed mesh and time step, the velocity matrix is independent of time if $l$ is constant. If $l(\phi)$ is evaluated explicitly, the matrix is known at each time step but may vary with time. 

\begin{remark}
In this paper, we treat the nonlinear bulk and boundary energies jointly through a unified scalar auxiliary variable. Instead of introducing a separate SAV for each contribution, this construction yields a simpler SAV–ZEC system with fewer auxiliary equations while retaining energy stability and a linear, decoupled solution structure. It also requires a lower bound only on the combined nonlinear energy, allowing greater flexibility than separate lower-bound assumptions on the two contributions.
\end{remark}

\section{Numerical experiments}
\label{sec:numerics}

This section is designed to verify the accuracy, stability, and applicability of the proposed method. Unless otherwise stated, the computational domain is $\Omega=[0,2]\times[-1,1]$, the generalized Navier and dynamic contact-line boundary conditions are imposed on the upper and lower walls, and periodic boundary conditions are used on the remaining boundaries for the phase-field, velocity,
and pressure. We also impose $\int_\Omega p\,\dd x=0$ to fix the additive constant in the pressure.

\subsection{Convergence test}

We first examine the temporal and spatial convergence behavior of the proposed method with stationary walls, namely, $\uw=0$.
The final time is chosen as $T=0.04$. 
Since no exact solution is available for this test, we use the adjacent-grid error to estimate the numerical convergence rate.

A sequence of uniformly refined meshes and time steps is used:
\[
h=\frac1{8},\frac1{16},\frac1{32},\frac1{64},\frac1{128},\frac1{256},
\qquad
\Delta t=\frac{0.04}{8},\frac{0.04}{16},\frac{0.04}{32},\frac{0.04}{64},\frac{0.04}{128},\frac{0.04}{256}.
\]
For each pair $(h,\Delta t)$, the numerical solution is computed up to $T=0.04$. 
The solution on the next finer mesh is taken as the reference solution after being restricted or interpolated onto the coarser mesh. 
For a generic variable $\xi$,  the adjacent-grid error in the
$L^2(\Omega)$ norm is defined by
\[
e_h(\xi)
=
\|\xi_h-\mathcal I_h \xi_{h/2}\|,
\]
where $\mathcal I_h$ denotes the interpolation or restriction operator from the finer grid to the coarser grid. 
The corresponding observed convergence order is computed by
\[
{\rm Order}
=
\log_2
\left(
\frac{e_h(\xi)}{e_{h/2}(\xi)}
\right).
\]
In the experiment, we compute the errors for the phase-field variable $\phi$, the velocity $\uh$, and the pressure $p$.

\begin{figure}[H]
	\centering
	\includegraphics[width=0.55\textwidth]{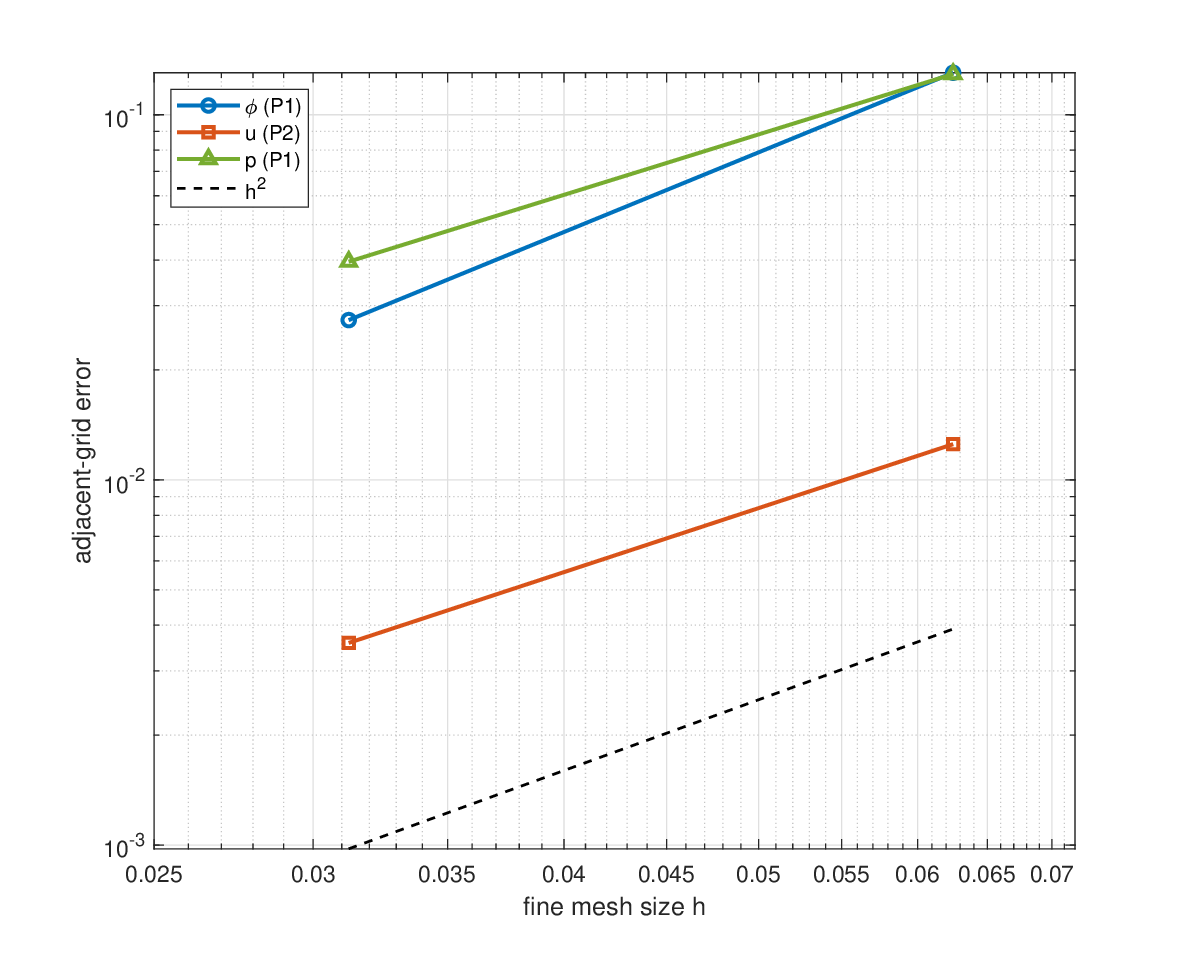}
	\caption{Adjacent-grid error convergence for $\phi$, $\uh$, and $p$ at $T=0.04$.}
	\label{fig:error-order}
\end{figure}
Figure~\ref{fig:error-order} shows approximately second-order
convergence for $\phi$ and $\uh$ under simultaneous space--time
refinement, consistent with second-order accuracy in time and
in the spatial $L^2$ norm.

\subsection{Two-phase Couette-flow benchmark}
\begin{figure}[H]
	\centering
	\includegraphics[width=0.9\textwidth]{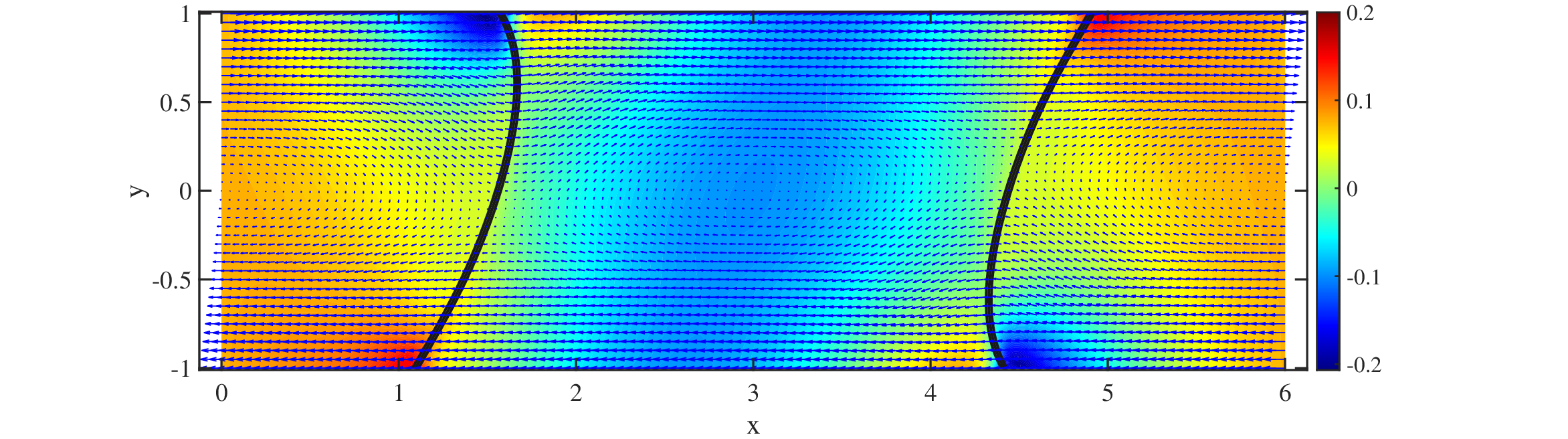}
	\caption{Two-phase Couette-flow benchmark: phase-field contour and velocity vectors at \(T=2\).}
	\label{fig:ex1-fields}
\end{figure}

We next consider a two-phase Couette-flow benchmark to illustrate the applicability of the proposed scheme under wall-driven shear. 
The computational domain is chosen as $\Omega=[0,6]\times[-1,1]$.
The initial velocity is taken as the standard Couette profile induced by the motion of the upper and lower walls. 
The initial phase-field function is prescribed by
\begin{equation*}
	\phi_0(x,y)
	=
	\tanh\left[
	\frac{1}{\sqrt{2}\varepsilon}
	\left(
	\frac{L_x}{4}
	-
	\left|x-\frac{L_x}{2}\right|
	\right)
	\right],
\end{equation*}
which describes a central strip bounded by two diffuse interfaces. 
The upper and lower walls move with opposite tangential velocities,
$\uw=(\pm 0.2,0)$, so that the flow is driven by shear.
The numerical parameters are chosen as
\[
\Delta t=2\times 10^{-4},
\qquad
h=\frac{1}{40},
\qquad
T=2.
\]
The physical and model parameters are
\[
M=0.01,\qquad
\lambda=1.2,\qquad
\varepsilon=0.05,\qquad
\gamma=100,\qquad
\theta_s=\frac{\pi}{3},
\]
and
\[
\nu=1,\qquad
\nu l(\phi)=\frac{1}{0.19},\qquad
A=1.
\]
Figure~\ref{fig:ex1-fields} shows the phase-field contour together with the velocity vectors at the final time \(T=2\). 
The two-phase interface is transported and deformed by the wall-driven shear flow, while the contact-line dynamics remains stable near the solid boundaries. 
The velocity field exhibits the expected Couette-type structure away from the interface, and local perturbations appear near the diffuse interfacial region due to the coupling between the phase field and the hydrodynamic motion.

\begin{figure}[H]
	\centering
	\begin{subfigure}{0.48\textwidth}
		\centering
		\includegraphics[width=\textwidth]{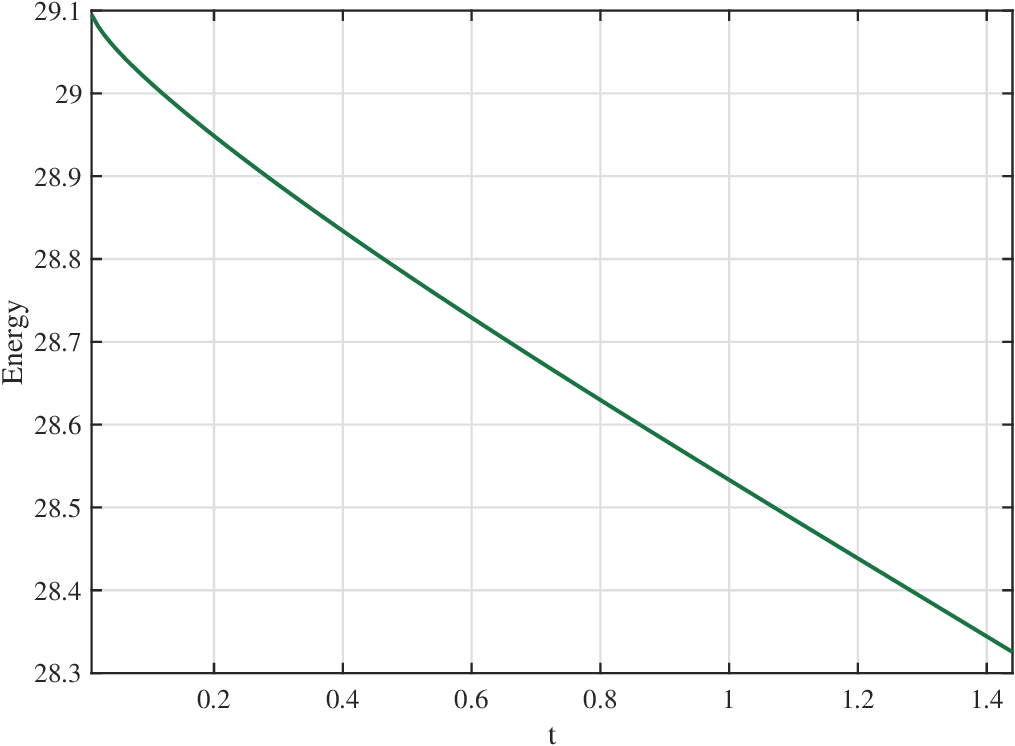}
		\caption{Modified discrete energy.}
	\end{subfigure}
	\hfill
	\begin{subfigure}{0.48\textwidth}
		\centering
		\includegraphics[width=\textwidth]{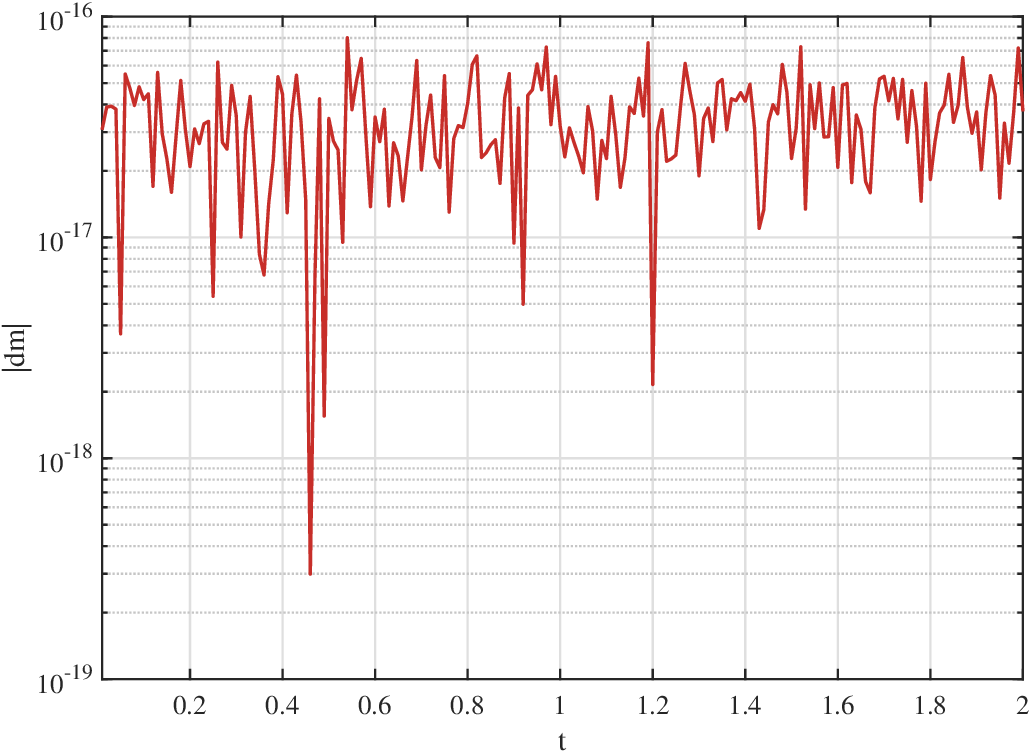}
		\caption{Mass variation.}
	\end{subfigure}
	\caption{Two-phase Couette-flow benchmark: energy behavior and mass conservation.}
	\label{fig:ex1-energy-mass}
\end{figure}

\begin{figure}[H]
	\centering
	\includegraphics[width=0.55\textwidth]{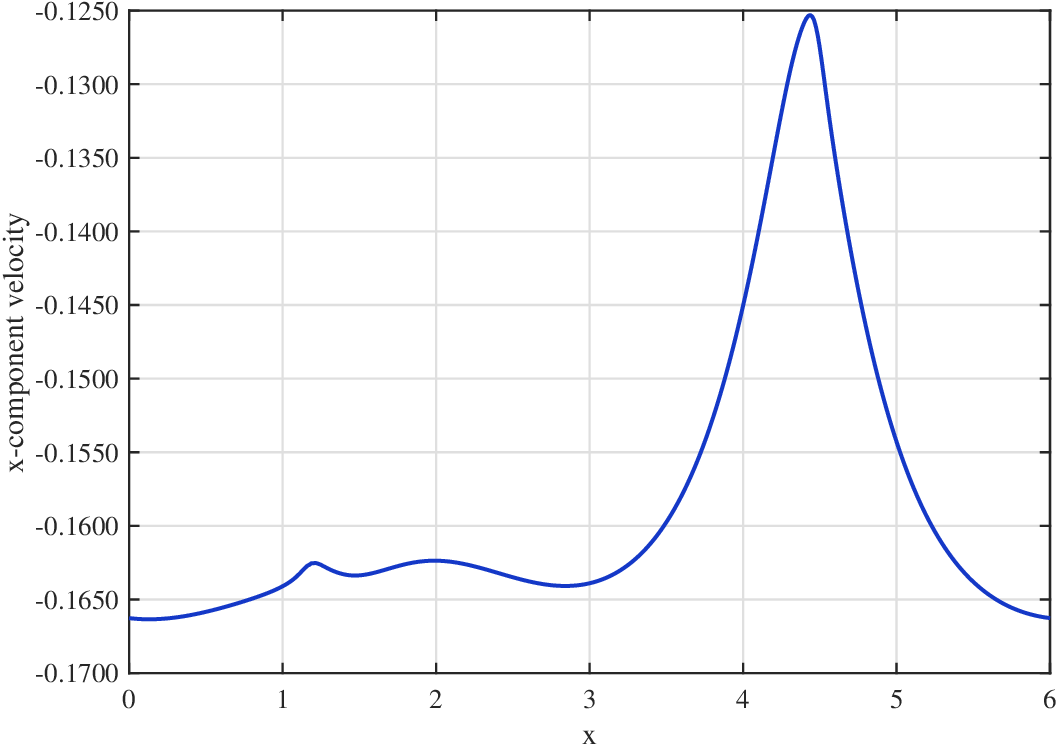}
	\caption{Two-phase Couette-flow benchmark: horizontal velocity profile along $y=-1$ at \(T=2\).}
	\label{fig:ex1-velocity-profile}
\end{figure}

We monitor the modified discrete energy and the total phase-field mass.
For the energy test, stationary walls, $\uw=0$, are imposed so that
no mechanical work enters through the boundary. As shown in
Figure~\ref{fig:ex1-energy-mass}, the modified energy decreases
monotonically, confirming the discrete energy-dissipation law. Moreover,
taking $\omega_h=1$ in \eqref{eq:scheme-phi} yields conservation of
$(\phi_h^{n+1},1)$. 
Accordingly, the absolute mass deviation $|(\phi_h^{n},1)-(\phi_h^{0},1)|$ remains at the level of roundoff error, consistent with the mass-conserving property of the phase-field update.
%
Finally, Figure~\ref{fig:ex1-velocity-profile} presents the horizontal velocity profile at \(T=2\). 
The numerical velocity follows the expected shear-driven profile induced by the opposite wall motions. 
The profile is nearly linear away from the diffuse interface, while deviations occur near the interfacial region because of capillary forcing and phase-field coupling. 

\subsection{Shear-driven droplet motion}
We next consider shear-driven droplet motion following
\cite{shen2015,yu2017}. The droplet is initially attached to the
lower wall, while the upper and lower walls move with opposite
tangential velocities. We use this test to examine the coupling
between wall slip, contact-line motion, and capillary forces
under shear. It therefore provides a representative benchmark for assessing the capability of the proposed fully decoupled method.

We set $\uh_0=0$ and prescribe the phase field by
\begin{equation*}
	\phi_0(x,y)
	=
	\tanh\left(
	\frac{\sqrt{(x-0.5L_x)^2+(y+1)^2}-1}
	{\sqrt{2}\varepsilon}
	\right),
\end{equation*}
which describes a diffuse droplet attached to the bottom boundary.
The computational parameters are chosen as
$\Delta t=10^{-3}, h_u=\frac{1}{40}, h_p=\frac{1}{20}$.
We use the Taylor–Hood \(P_2/P_1\) finite element pair for the velocity and pressure. 
The model parameters are $M=0.001,\lambda=12,
\varepsilon=0.05, \gamma=100, \nu=1, \nu l(\phi)=\frac{1}{0.19}.	$

\begin{figure}[htbp]
	\centering
	\setlength{\tabcolsep}{0pt}
	\renewcommand{\arraystretch}{0}
	
	\begin{tabular}{@{}c@{\hspace{0.8mm}}c@{\hspace{0.8mm}}c@{}}
		\includegraphics[
		width=0.33\textwidth,
		trim={25bp 18bp 6bp 4bp},
		clip
		]{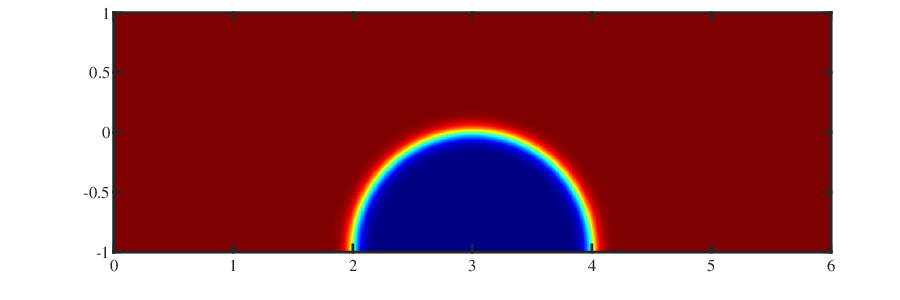}
		&
		\includegraphics[
		width=0.33\textwidth,
		trim={25bp 18bp 6bp 4bp},
		clip
		]{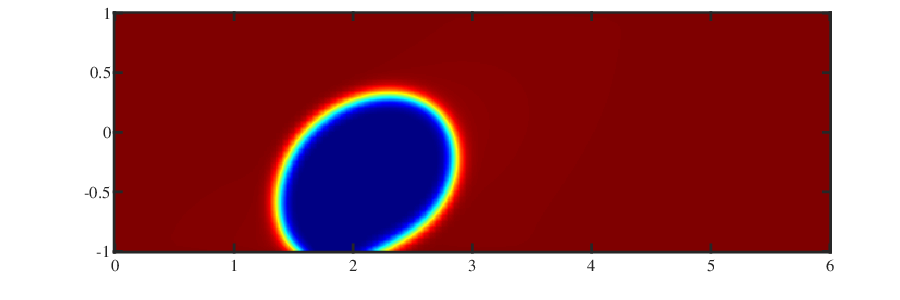}
		&
		\includegraphics[
		width=0.33\textwidth,
		trim={25bp 18bp 6bp 4bp},
		clip
		]{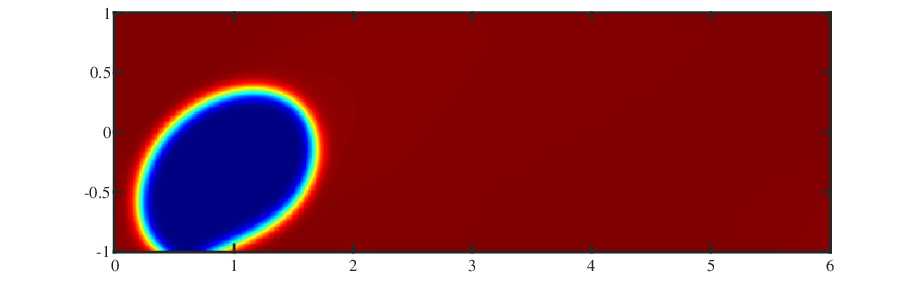}
		\\[0.8mm]
		\includegraphics[
		width=0.33\textwidth,
		trim={25bp 18bp 6bp 4bp},
		clip
		]{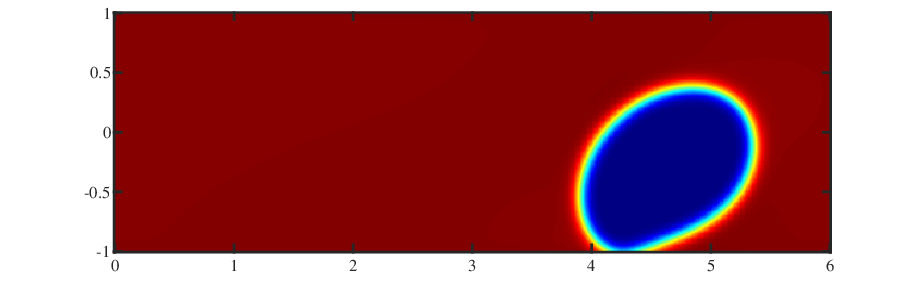}
		&
		\includegraphics[
		width=0.33\textwidth,
		trim={25bp 18bp 6bp 4bp},
		clip
		]{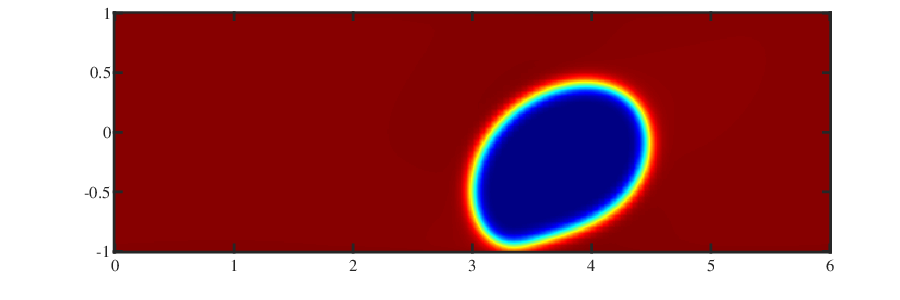}
		&
		\includegraphics[
		width=0.33\textwidth,
		trim={25bp 18bp 6bp 4bp},
		clip
		]{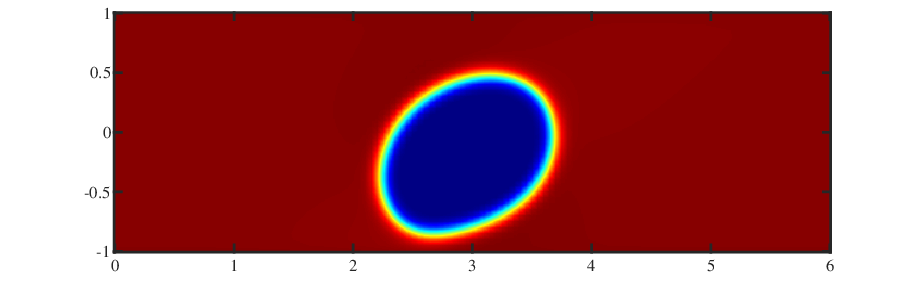}
	\end{tabular}
	
	\vspace{-2mm}
	\caption{Shear-driven droplet motion for $\theta_s=30^\circ$ at
		$t=0,\,1,\,2,\,4,\,5,$ and $6$, ordered from left to right
		and top to bottom.}
	\label{fig:shear_droplet_30}
\end{figure}

\begin{figure}[htbp]
	\centering
	\setlength{\tabcolsep}{0pt}
	\renewcommand{\arraystretch}{0}
	
	\begin{tabular}{@{}c@{\hspace{0.8mm}}c@{\hspace{0.8mm}}c@{}}
		\includegraphics[
		width=0.33\textwidth,
		trim={25bp 18bp 6bp 4bp},
		clip
		]{t=0.eps}
		&
		\includegraphics[
		width=0.33\textwidth,
		trim={25bp 18bp 6bp 4bp},
		clip
		]{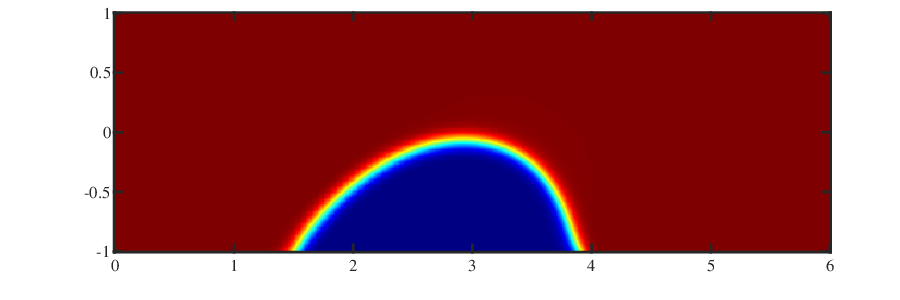}
		&
		\includegraphics[
		width=0.33\textwidth,
		trim={25bp 18bp 6bp 4bp},
		clip
		]{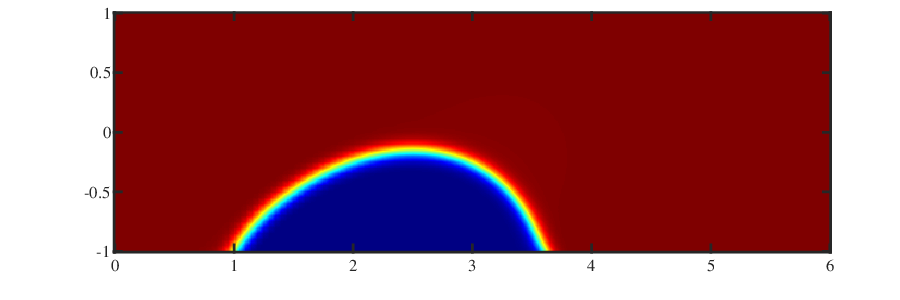}
		\\[0.8mm]
		\includegraphics[
		width=0.33\textwidth,
		trim={25bp 18bp 6bp 4bp},
		clip
		]{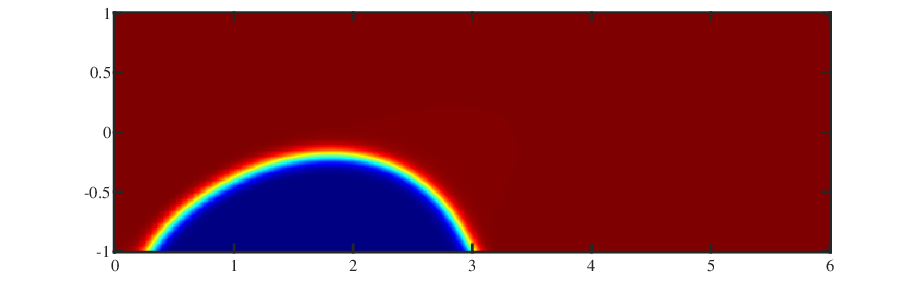}
		&
		\includegraphics[
		width=0.33\textwidth,
		trim={25bp 18bp 6bp 4bp},
		clip
		]{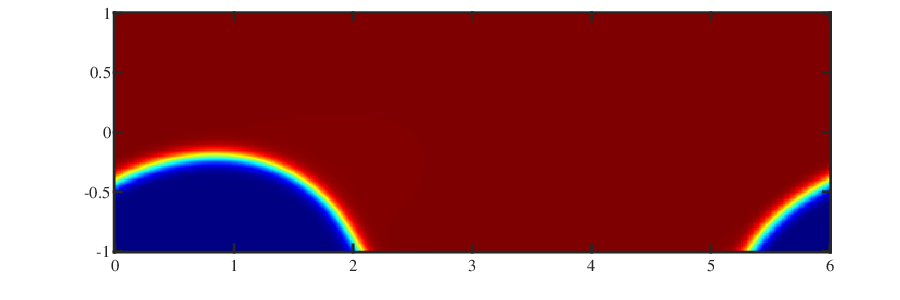}
		&
		\includegraphics[
		width=0.33\textwidth,
		trim={25bp 18bp 6bp 4bp},
		clip
		]{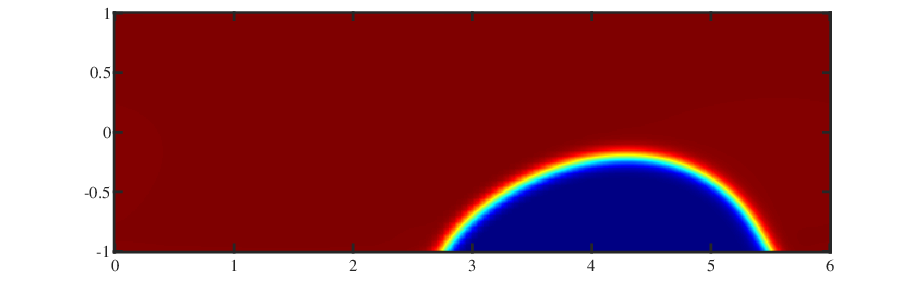}
	\end{tabular}
	
	\vspace{-2mm}
	\caption{Shear-driven droplet motion for $\theta_s=120^\circ$ at
		$t=0,\,0.3,\,0.6,\,1,\,1.5,$ and $3$, ordered from left to right
		and top to bottom.}
	\label{fig:shear_droplet_120}
\end{figure}

We first consider the case $\theta_s=30^\circ$.
As shown in Figure~\ref{fig:shear_droplet_30}, the initially
wall-attached droplet becomes tilted and elongated under the imposed shear. It is subsequently lifted away from the lower wall and evolves into a detached droplet in the bulk flow. For $\theta_s=120^\circ$, Figure~\ref{fig:shear_droplet_120} shows a
distinctly different phenomenon. The droplet remains attached to the lower wall, its contact region widens and  moves along the substrate. Its passage through the periodic boundary is visible at intermediate times, while no apparent detachment occurs.

\subsection{Three-dimensional shear-driven droplet dynamics}
To further assess the performance of the proposed scheme in three dimensions, we consider a three-dimensional shear-driven droplet test.  The computational domain is $\Omega=(0,3)\times(0,1)^2$.
The initial droplet is a hemisphere of radius \(R_0=0.4\), and is placed on the lower substrate. 
The initial velocity is set to zero, and the phase-field distribution is initialized by
\[
\phi_0(x,y,z)
=
\tanh\left(
\frac{\sqrt{(x-2)^2+(y-0.5)^2+z^2}-0.4}
{\sqrt{2}\varepsilon}
\right).
\]
The droplet is driven by opposite motions of the upper and lower plates,
\[
\uw=(0.5,0,0)^T \,\text{on}\, z=1,\qquad
\uw=(-0.5,0,0)^T \,\text{on}\, z=0.
\]
The mesh size and time step are chosen as $
h=0.05$, $ \Delta t=0.0025$. 

\begin{figure}[H]
	\centering
	\includegraphics[
	width=1.05\textwidth,
	trim={55bp 35bp 45bp 55bp},
	clip
	]{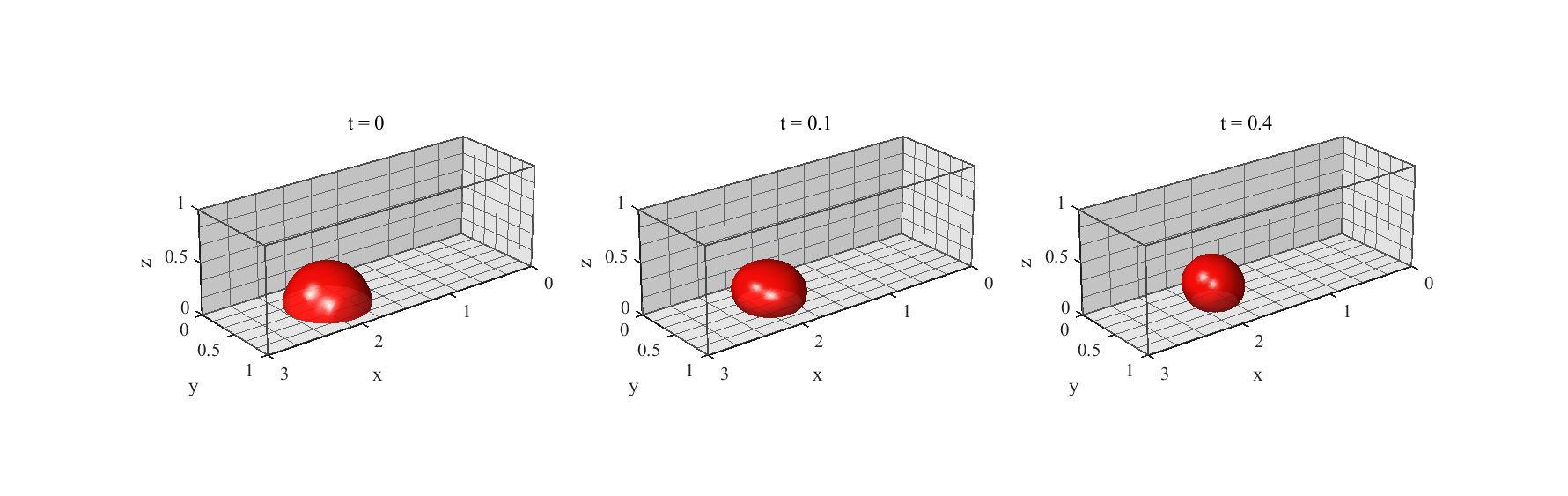}
	\vspace{-3mm}
	\caption{Three-dimensional shear-driven droplet dynamics with
		$\theta_s=30^\circ$. Snapshots at $t=0,\,0.1,$ and $0.4$
		show the deformation and migration of the droplet under the
		imposed shear flow.}
	\label{fig:3d-shear-30}
\end{figure}

We first set the equilibrium contact angle to
$\theta_s=30^\circ$. As shown in
Figure~\ref{fig:3d-shear-30}, the initially hemispherical droplet is
deformed and displaced in the streamwise direction by the imposed
shear. At $t=0.1$, it becomes horizontally elongated and asymmetric.
By $t=0.4$, its contact region on the lower substrate has contracted,
while the droplet becomes taller and approaches a nearly spherical
shape. The interface remains smooth during this deformation and
migration process.

\begin{figure}[H]
	\centering
	\includegraphics[
	width=1.05\textwidth,
	trim={55bp 220bp 45bp 200bp},
	clip
	]{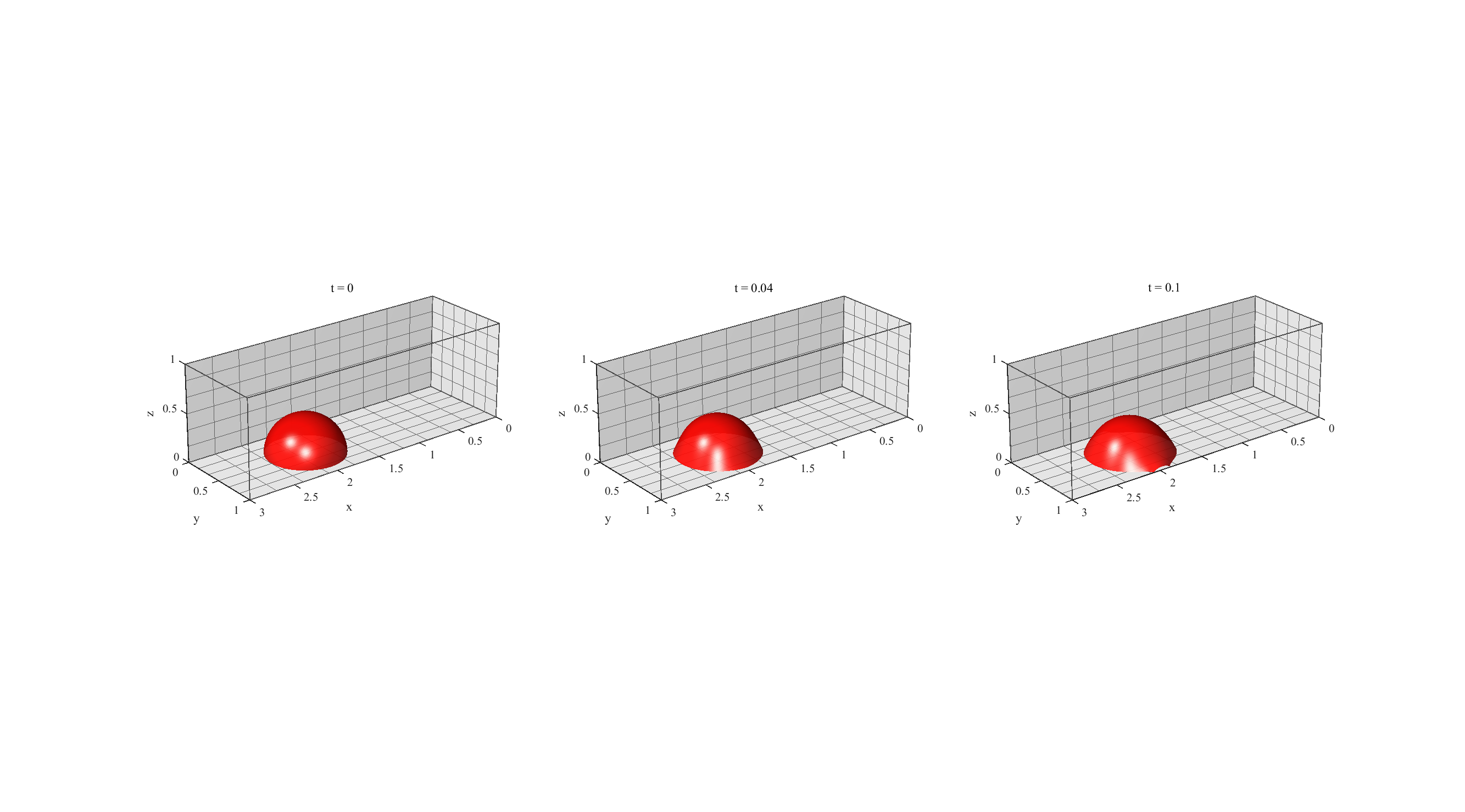}
	\vspace{-2mm}
	\caption{Three-dimensional shear-driven droplet dynamics with
		$\theta_s=120^\circ$. Snapshots at $t=0,\,0.04,$ and $0.1$
		show the deformation, migration, and progressive increase in
		the droplet--substrate contact area under the imposed shear flow.}
	\label{fig:3d-shear-120}
\end{figure}

We then repeat the simulation with $\theta_s=120^\circ$.
Figure~\ref{fig:3d-shear-120} shows a markedly different evolution.
The droplet remains attached to the lower substrate and progressively spreads in the streamwise direction. Its contact area increases and its height decreases over $0\leq t\leq 0.1$, resulting in a wider and flatter cap-shaped configuration. These results suggest that the proposed scheme can capture
differences in three-dimensional contact-line motion for
different prescribed contact angles.

\subsection{Relaxation on a curved surface}
Finally, we consider droplet relaxation with stationary walls
in a non-rectangular domain whose curved lower wall is approximated
by a polygonal boundary. This test examines droplet evolution
and discrete energy dissipation in this geometry.

The computational domain is given by
\[
-0.5\le x\le 0.5,\qquad
\max\{-0.2,-1.65+\sqrt{1.5^2-x^2}\}\le y\le 0.2.
\]
The left and right boundaries are treated periodically. 
We impose the dynamic contact-line condition on the upper and lower walls. 
The wall velocity is set to zero, \(\uw=0\). 
Therefore, the boundary term involving $\uw$ vanishes, and the modified discrete energy is expected to decay monotonically. The parameters are chosen as
\[
M=5\times 10^{-4},\quad 
\lambda=2.4,\quad 
\varepsilon=0.01,\quad
\gamma=200,\quad 
\nu=0.2,\quad 
l(\phi)=80,\quad
\theta_s=\frac{2\pi}{3}.
\]
The time step and final time are $\Delta t=10^{-4}$, $T=1$. The initial phase field is given by
\[
\phi_0(x,y)=
\tanh\left(
\frac{0.2-\sqrt{x^2+(y+0.15)^2}}{\sqrt{2}\varepsilon}
\right),
\]
which describes a circular droplet in contact with the curved
lower wall. The droplet is initially out of equilibrium and
subsequently relaxes under capillary forces and the dynamic
contact-line condition.

\begin{figure}[htbp]
	\centering
	\setlength{\tabcolsep}{0pt}
	\renewcommand{\arraystretch}{0}
	
	\begin{tabular}{@{}c@{\hspace{0.8mm}}c@{}}
		\includegraphics[
		width=0.48\textwidth,
		trim={25bp 18bp 6bp 4bp},
		clip
		]{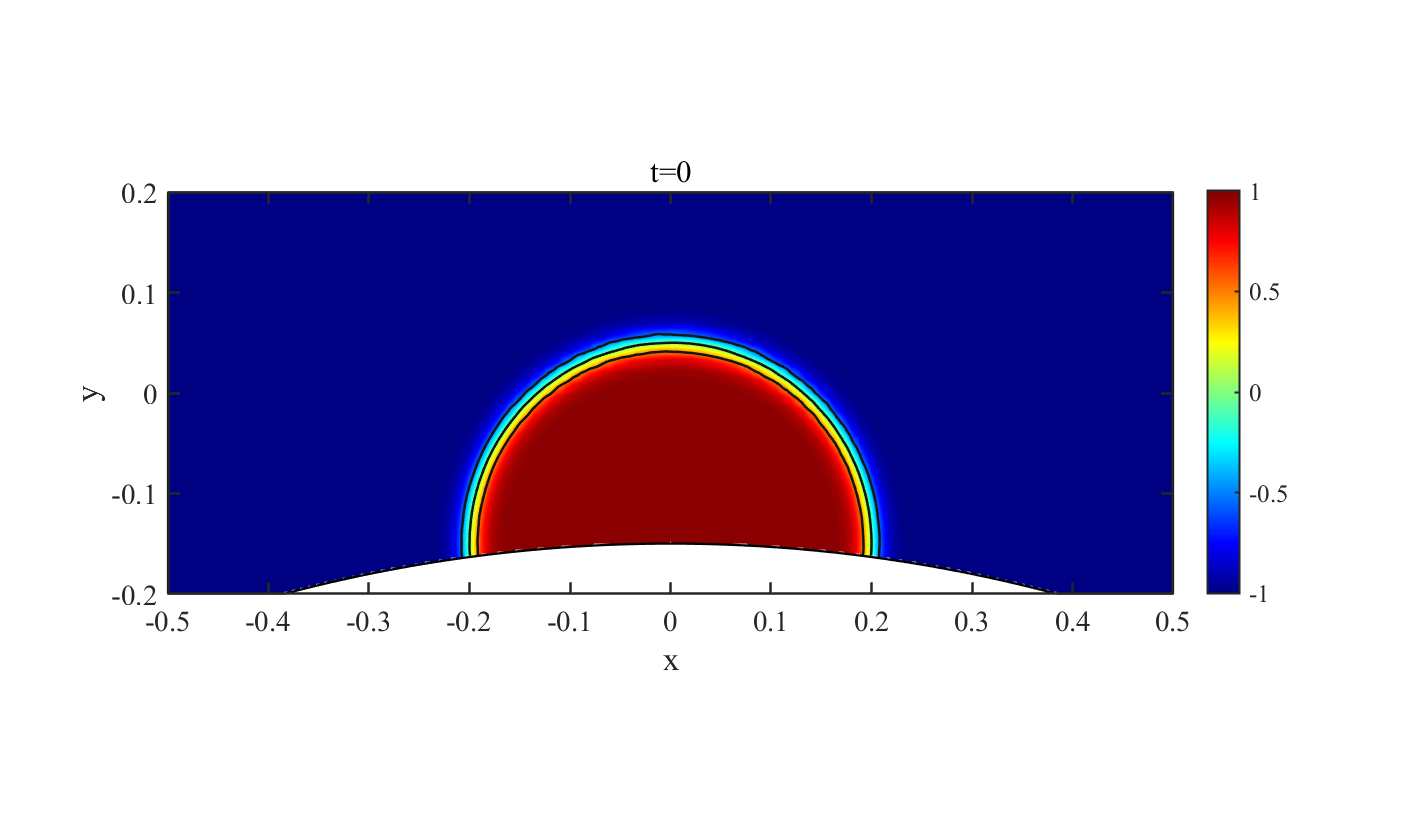}
		&
		\includegraphics[
		width=0.48\textwidth,
		trim={25bp 18bp 6bp 4bp},
		clip
		]{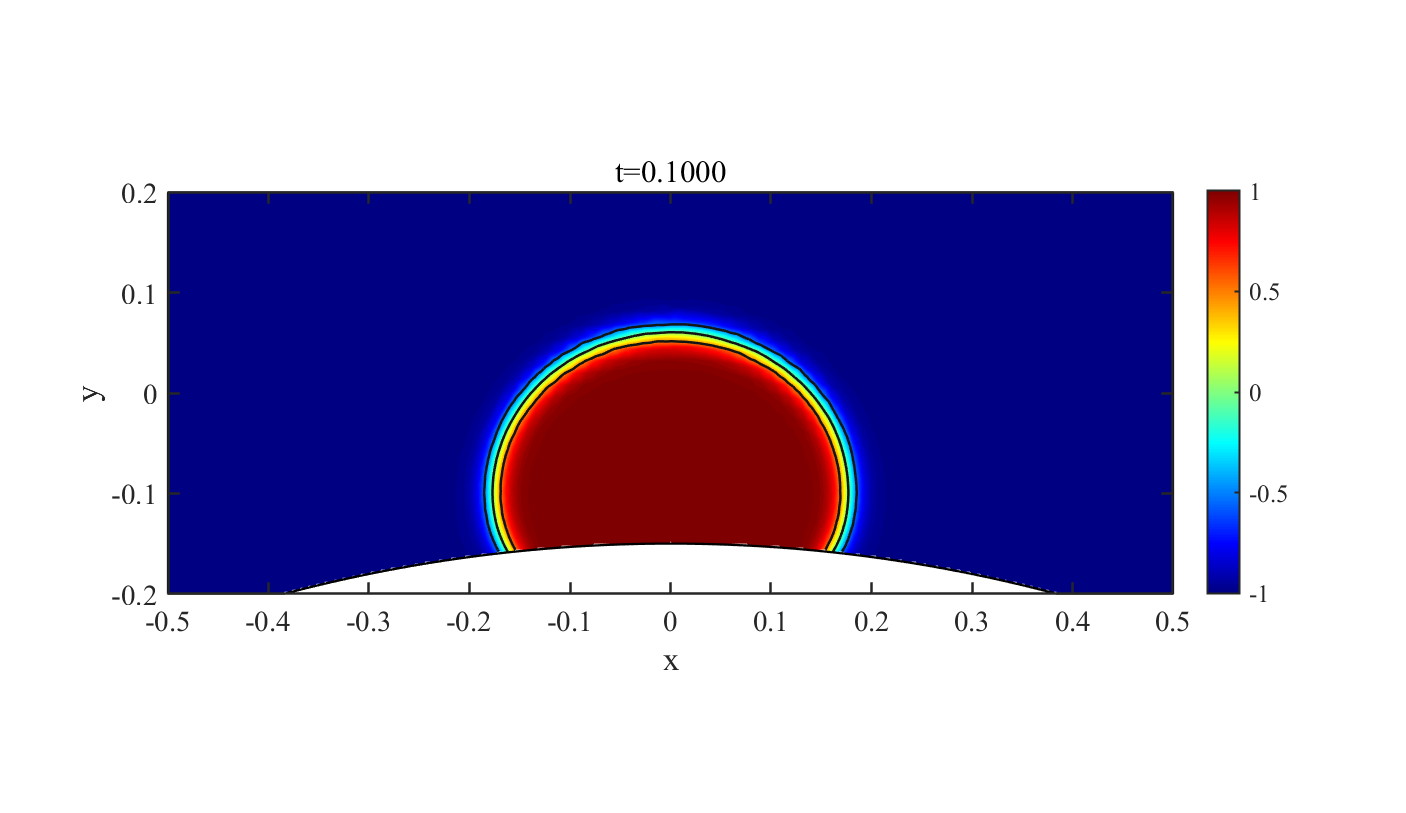}
		\\[0.8mm]
		\includegraphics[
		width=0.48\textwidth,
		trim={25bp 18bp 6bp 4bp},
		clip
		]{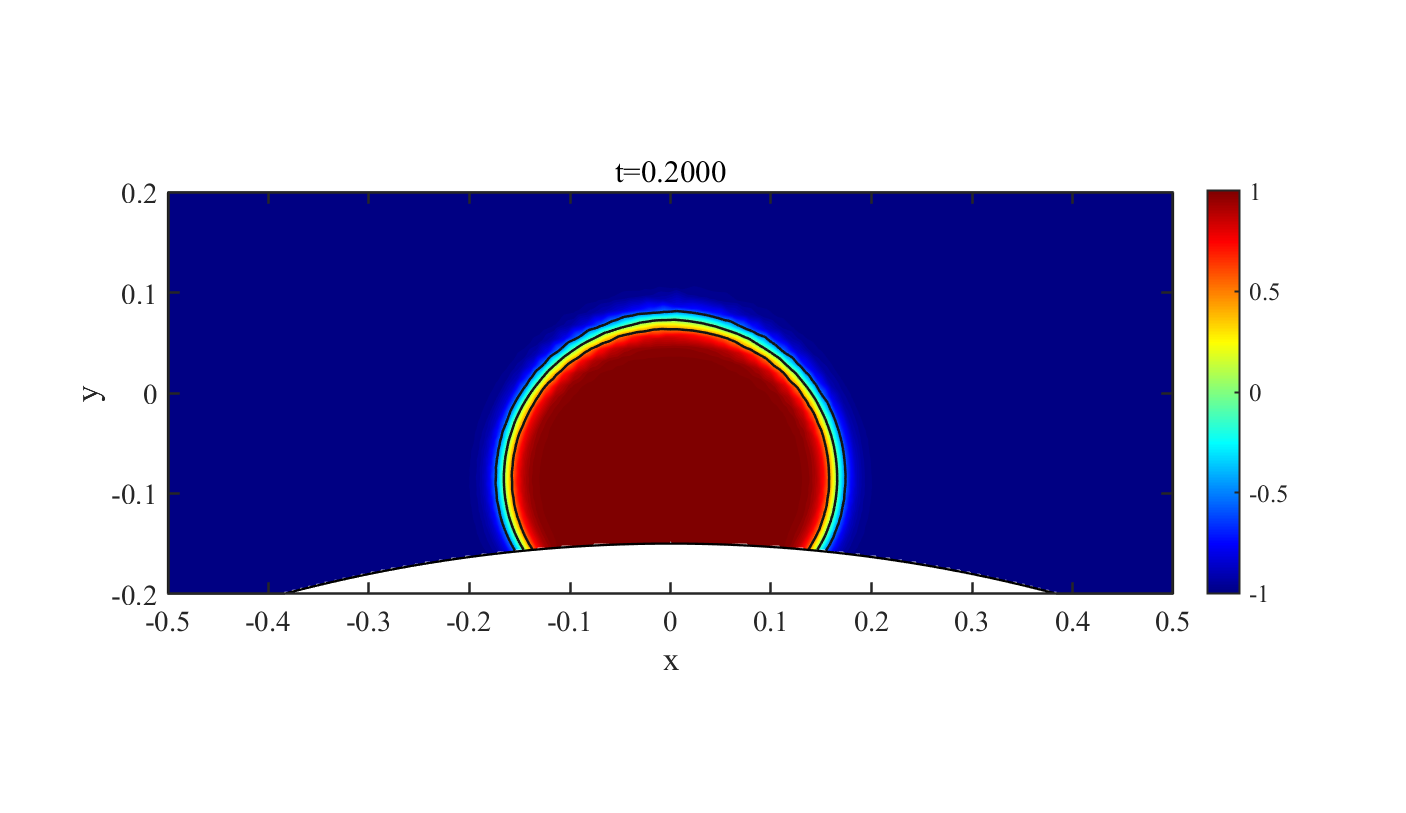}
		&
		\includegraphics[
		width=0.48\textwidth,
		trim={25bp 18bp 6bp 4bp},
		clip
		]{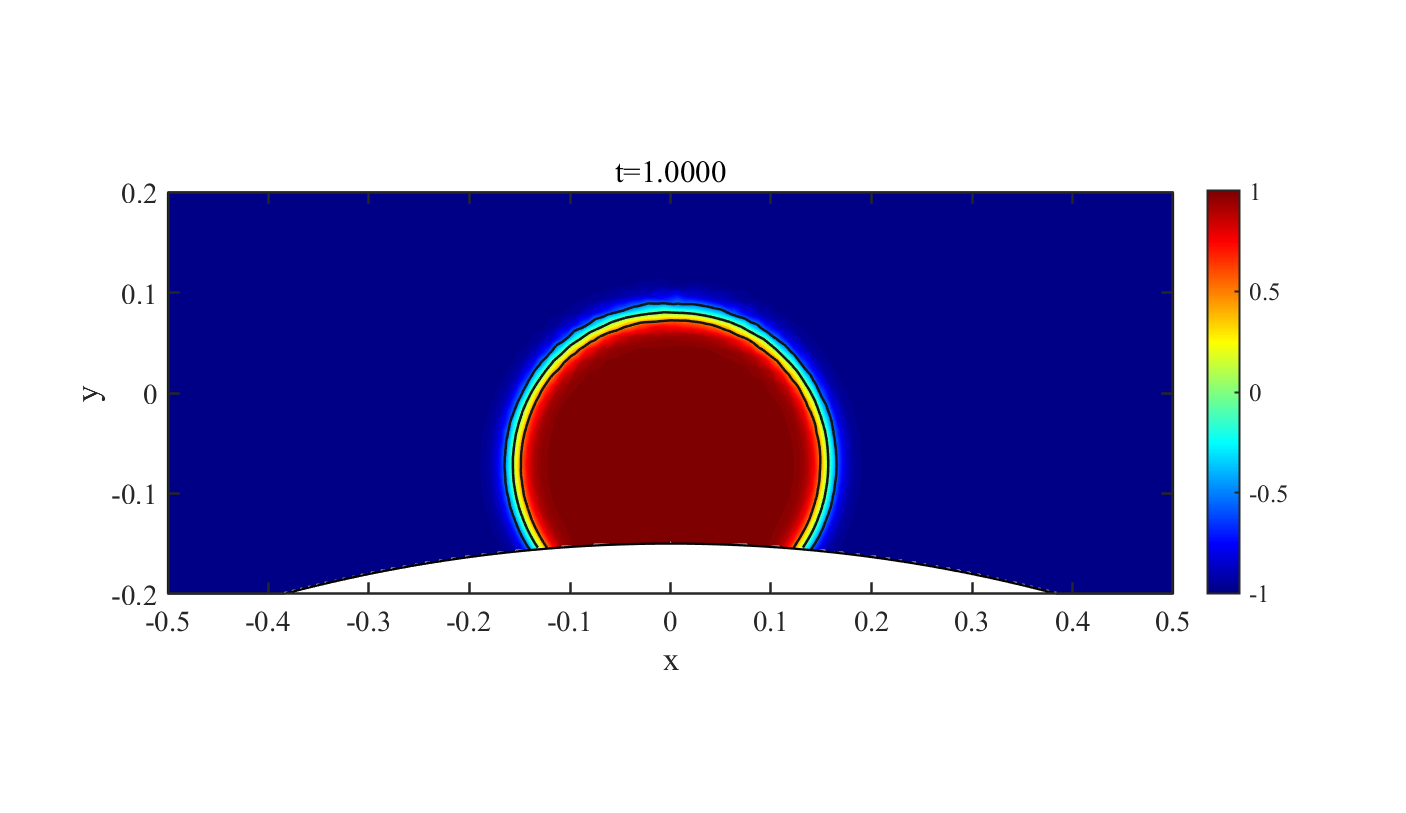}
	\end{tabular}
	
	\vspace{-2mm}
	\caption{Relaxation on a polygonal approximation of a curved surface:
		phase-field snapshots at \(t=0,0.1,0.2,\) and \(1\).
		The black curve denotes the zero level set of \(\phi\).}
	\label{fig:ex2-polygon-snapshots}
\end{figure}

Figure~\ref{fig:ex2-polygon-snapshots} shows the phase field at \(t=0,0.1,0.2,\) and \(1\). The black curve is the zero level set of \(\phi\),
making the transient deformation and oscillatory adjustment induced by the
curved substrate more clearly visible. Starting from the initial circular
shape, the droplet rapidly adjusts its contact region on the curved wall.
During the early stage, the contact line moves along the lower boundary,
while the interface undergoes noticeable deformation and a damped
oscillatory response. As time increases, these oscillations gradually
decay, the contact-line motion slows down and the interface changes only slightly. This
behavior is consistent with the contact-line relaxation mechanism described
by the DCLC on a curved substrate.

\begin{figure}[htbp]
	\centering
	\setlength{\tabcolsep}{0pt}
	\renewcommand{\arraystretch}{0}
	
	\begin{tabular}{@{}c@{\hspace{0.8mm}}c@{}}
		\includegraphics[
		width=0.48\textwidth,
		trim={25bp 18bp 6bp 4bp},
		clip
		]{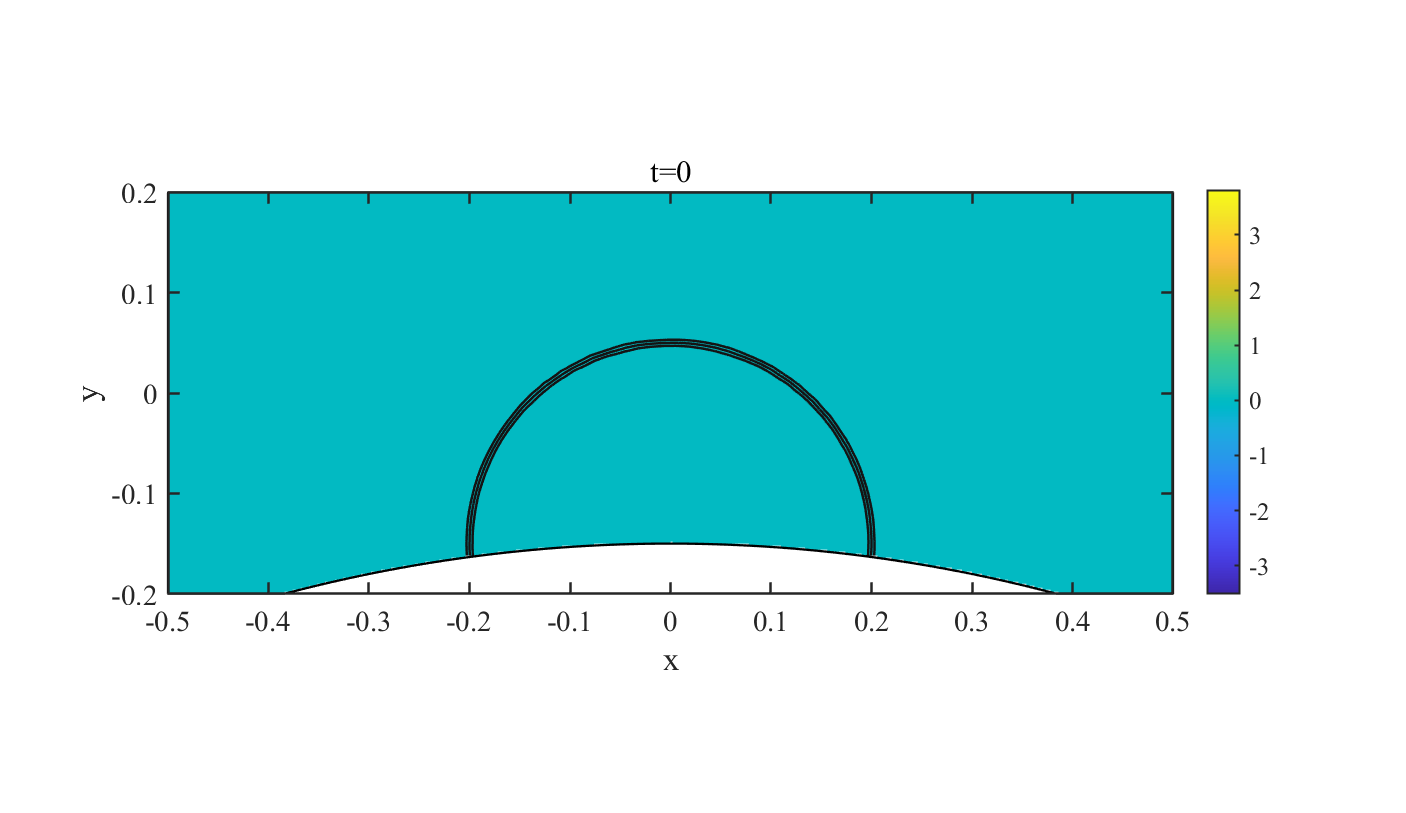}
		&
		\includegraphics[
		width=0.48\textwidth,
		trim={25bp 18bp 6bp 4bp},
		clip
		]{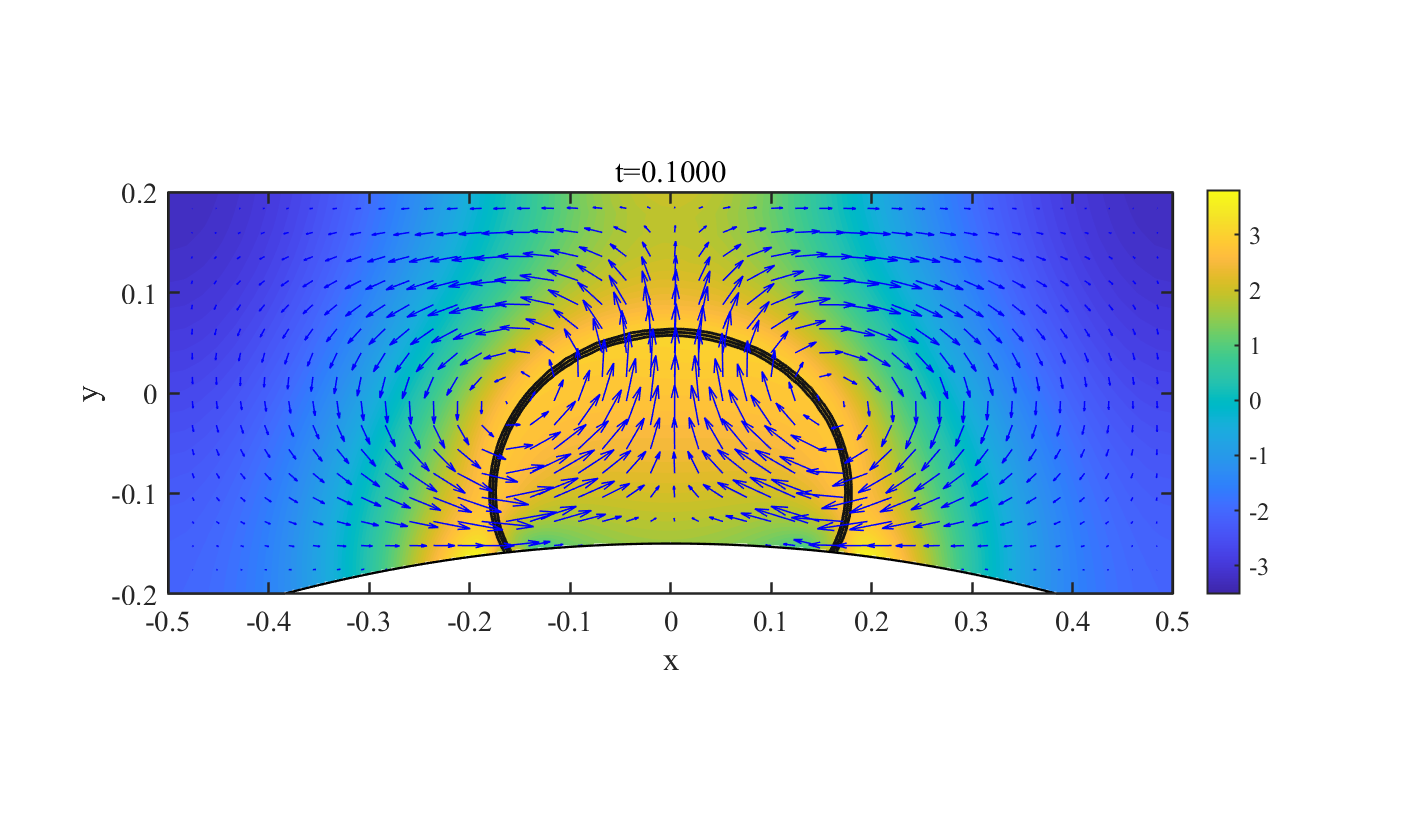}
		\\[-0.5mm]
		\includegraphics[
		width=0.48\textwidth,
		trim={25bp 18bp 6bp 4bp},
		clip
		]{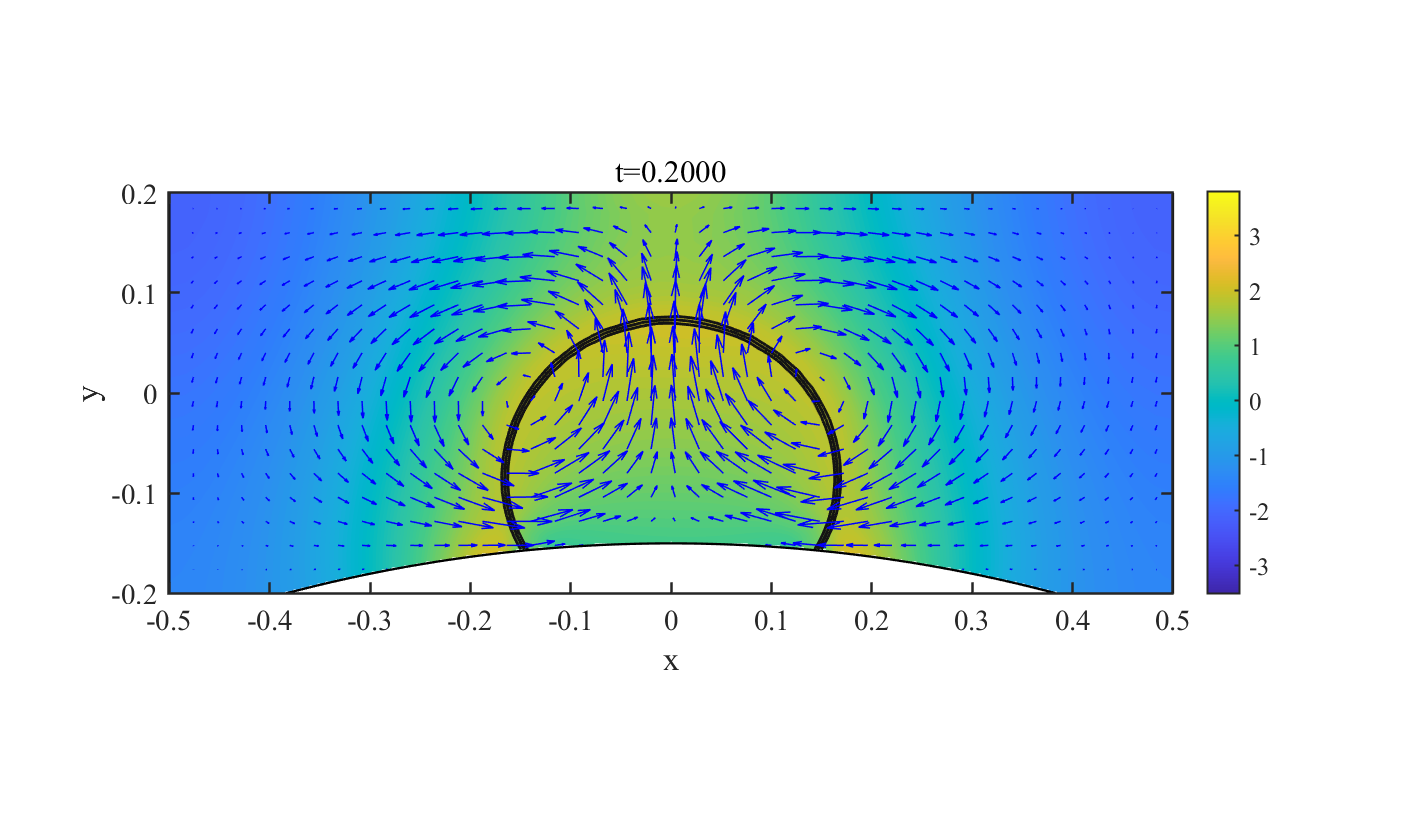}
		&
		\includegraphics[
		width=0.48\textwidth,
		trim={25bp 18bp 6bp 4bp},
		clip
		]{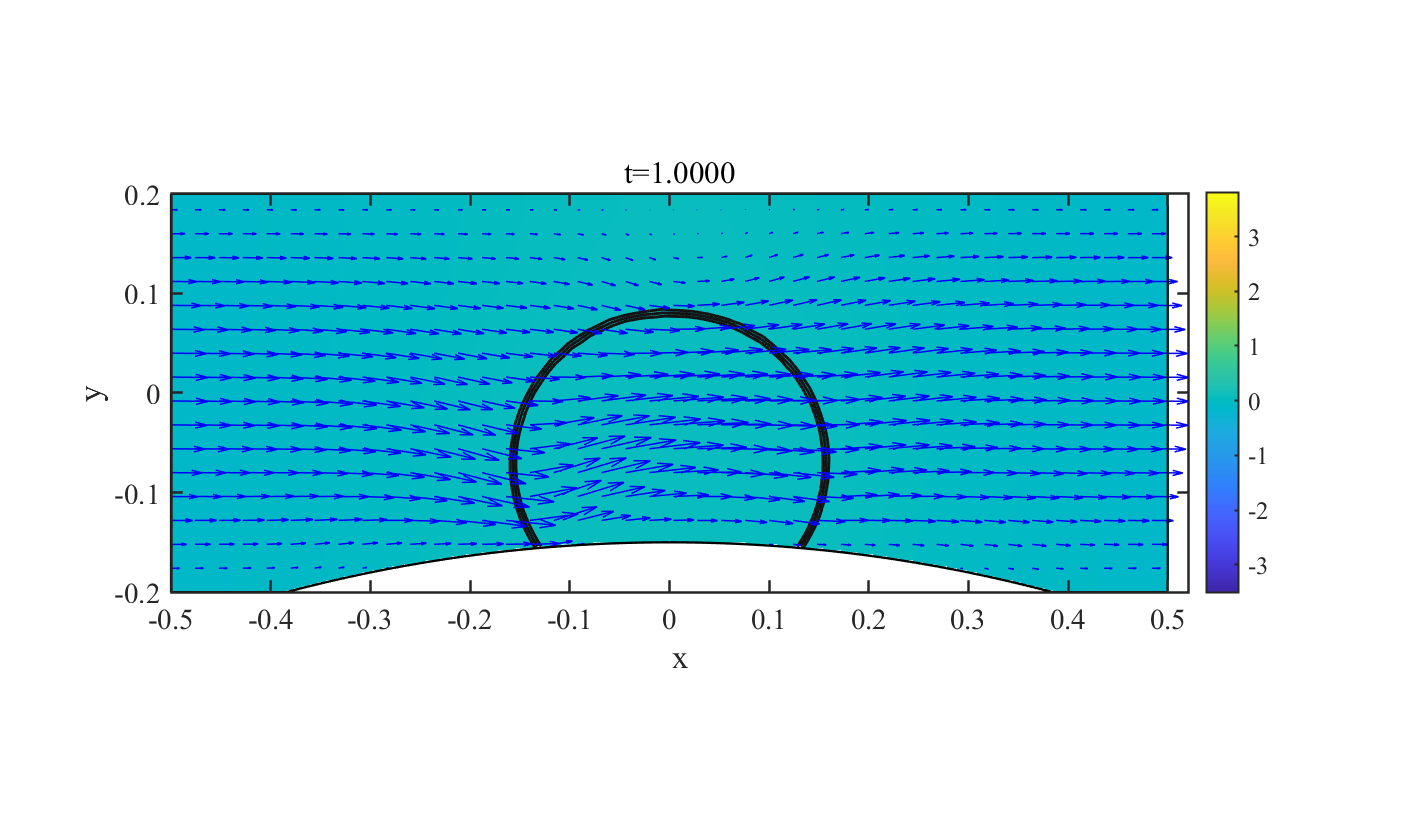}
	\end{tabular}
	
	\vspace{-2mm}
	\caption{Relaxation on a polygonal approximation of a curved surface:
		Pressure, phase-field contours, and velocity  at \(t=0,0.1,0.2,\) and \(1\).
		The color denotes the pressure \(p\), the black curves denote selected
		phase-field contours, and the blue arrows denote the velocity field.}
	\label{fig:ex2-polygon-coupled}
\end{figure}

Figure~\ref{fig:ex2-polygon-coupled} displays the corresponding pressure field, phase-field contours, and velocity vectors. 
At early times, pressure variations are concentrated near the interface and the contact region. 
The induced velocity field is also stronger during this initial relaxation stage. 
As the droplet approaches its equilibrium shape, both the pressure variation and the velocity magnitude decrease significantly. 
At \(t=1\), the velocity field is nearly stationary, indicating that the droplet has almost reached a relaxed state on the curved wall.

\begin{figure}[H]
	\centering
	\includegraphics[width=0.55\textwidth]{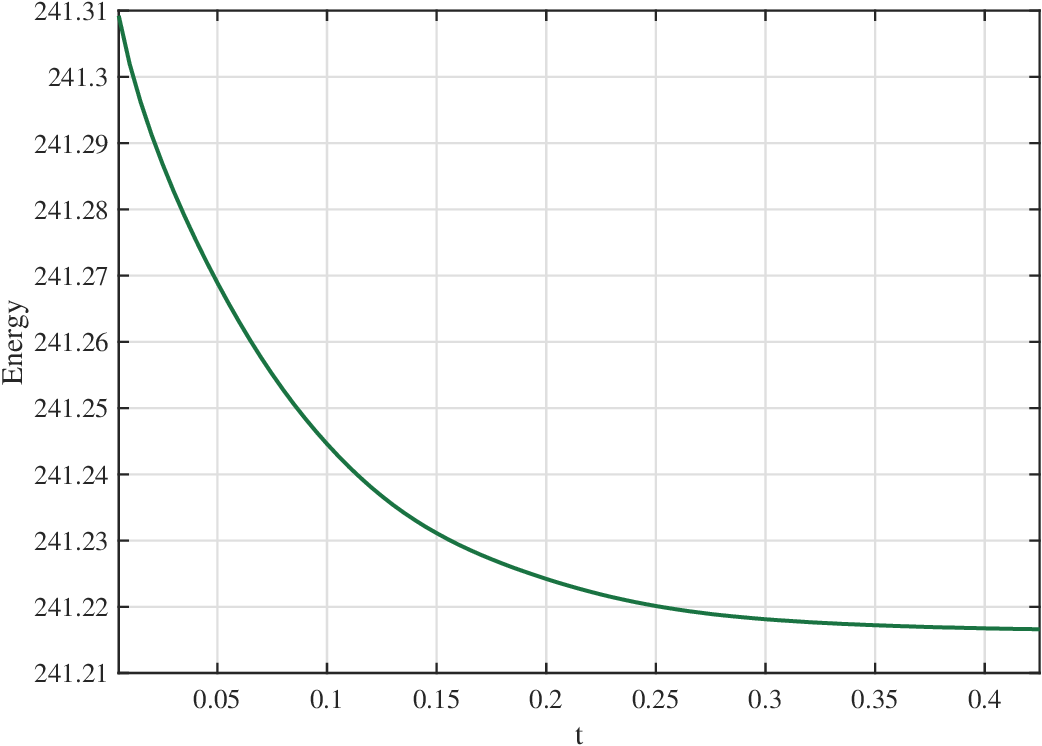}
	\caption{Relaxation on a polygonal approximation of a curved surface: monotone decay of the modified discrete energy.}
	\label{fig:ex2-polygon-energy}
\end{figure}

Since the wall is static in this test, the theoretical energy law predicts monotone energy dissipation. As shown in Figure~\ref{fig:ex2-polygon-energy}, the modified discrete energy decreases
monotonically from $241.3092$ to $241.2161$, in agreement with
Theorem~\ref{Thm5.1}.

\section{Conclusion}
\label{sec:conclusion}
We have developed a linear, fully decoupled BDF2 scheme for the CH--NS moving-contact-line model. The scheme introduces a unified SAV for the nonlinear bulk and boundary energies and a ZEC
variable for the nonlinear coupling terms. This unified SAV–ZEC system reduces the number of auxiliary equations while preserving energy stability and a linear, decoupled solution structure. This formulation also requires only the combined nonlinear energy to be bounded from below.
%
We have proved unconditional discrete energy stability for stationary walls.  Numerical tests show approximately second-order
convergence under simultaneous space--time refinement. 
Further simulations illustrate shear-driven droplet motion and relaxation on flat and curved walls, demonstrating the applicability of the finite element scheme to non-rectangular domains.
Future work includes error estimates, adaptive mesh refinement,
and extensions to variable densities and viscosities.

\section*{Acknowledgments}
The authors thank Professor Xiaofeng Yang of the University
of South Carolina for helpful discussions and valuable
suggestions on this work. X. Bao is partially supported by the NSFC (No. 12201050) and by Science and
Technology Projects in Guangzhou (No. SL2024A04J00931).

\bibliographystyle{elsarticle-num}
\bibliography{references}
\end{document}